\documentclass[a4paper,12pt,headsepline,bibliography=totoc]{scrartcl}
\usepackage[latin1]{inputenc}
\usepackage[T1]{fontenc}
\usepackage{lmodern}
\usepackage[english]{babel}
\usepackage{amsthm}
\usepackage{amsfonts}  
\usepackage{amsrefs}
\usepackage{cite} 
\usepackage{amsmath}
\usepackage{amssymb} 
\usepackage{mathtools}
\usepackage{lipsum}
\usepackage{nicefrac}
\usepackage{pdfpages} 
\usepackage{titleref}
\usepackage{paralist}
\usepackage{ngerman}
\usepackage{graphicx}
\usepackage{verbatim}
\usepackage{amsbsy}
\usepackage{fancyhdr}
\usepackage{bbm}
\usepackage{enumitem}
\usepackage{tabto}
\usepackage{bm}
\usepackage{array}
\usepackage{booktabs}
\usepackage{multirow}
\usepackage{bigdelim} 
\usepackage{rkeyval}
\usepackage{authblk}

\newtheoremstyle{mystyle}
  {.35cm}% measure of space to leave above the theorem. E.g.: 3pt
  {.35cm}% measure of space to leave below the theorem. E.g.: 3pt
  { \nopagebreak\itshape}% name of font to use in the body of the theorem
  {}% measure of space to indent
  {\bfseries}% name of head font
  {}% punctuation between head and body
  {\newline}% space after theorem head; " " = normal interword space
  {}% Manually specify head

\theoremstyle{mystyle}

\newtheorem{T}{Theorem}[section]
\newtheorem{D}[T]{Definition}

\newtheorem{N}[T]{Notation}
\newtheorem{Rem}[T]{Remark}
\newtheorem{Prop}[T]{Proposition}
\newtheorem{Le}[T]{Lemma}
\newtheorem{Co}[T]{Corollary}
\newtheorem{C}[T]{Assumption}

\providecommand{\Rd}[1]{\mathbb{R}^{#1}}
\providecommand{\norm}[1]{\left\lVert#1\right\rVert}
\providecommand{\abs}[1]{\left\lvert#1\right\rvert}

\providecommand{\sg}[2]{\left(T_{t#2}^{#1}\right)_{t\geq0}}
\providecommand{\asg}[2]{\left(\hat{T}_{t#2}^{#1}\right)_{t\geq0}}

\newcounter{subequation}
\newlength\mtabskip\mtabskip=-1.25cm

\def\mtabLong{long}


\numberwithin{equation}{section}

\title{Overdamped limit of generalized stochastic Hamiltonian systems for singular interaction potentials}

\author[1]{Martin Grothaus}
\author[2]{Andreas Nonnenmacher}
\affil[1]{Mathematics Department, TU Kaiserslautern,\newline PO Box 3049, 67653 Kaiserslautern, Germany (grothaus@mathematik.uni-kl.de)}
\affil[2]{Mathematics Department, TU Kaiserslautern,\newline PO Box 3049, 67653 Kaiserslautern, Germany (nonnenmacher@mathematik.uni-kl.de)}
\begin{document}
\newcounter{PHIeins}
\stepcounter{PHIeins}
\newcounter{PHIzwei}
\stepcounter{PHIzwei}

%\title{Essential M-Dissipativity of the Generator of
%	a generalized Stochastic Hamiltonian System and 
%	its scaling limit}

\maketitle
\begin{abstract}
First weak solutions of generalized stochastic Hamiltonian systems (gsHs) are constructed via essential m-dissipativity of their generators on a suitable core. For a scaled gsHs we prove convergence of the corresponding semigroups and tightness of the weak solutions. This yields convergence in law of the scaled gsHs to a distorted Brownian motion. In particular, the results confirm the convergence of the Langevin dynamics in the overdamped regime to the overdamped Langevin equation.
The proofs work for a large class of (singular) interaction potentials including, e.g., potentials of Lennard--Jones type.


\end{abstract}

\thispagestyle{empty}

% !TeX spellcheck = en_US
\section[Introduction]{Introduction}

The motion of interacting particles in a surrounding medium can be described by the Langevin equation, i.e.,
\begin{align}
dX_t&=V_tdt,\tag{1.1a}\label{Langevin1}\\
dV_t&=-\nabla\Phi_1(X_t)dt-\gamma V_tdt+\sqrt{2\gamma\beta^{-1}}dB_t,\tag{1.1b}\label{Langevin2}
\intertext{where $ \nabla \Phi_1 $ prescribes external and interacting forces between the particles, $ \gamma>0 $ is a constant describing the magnitude of friction, $ \beta>0 $ is up to a constant the inverse temperature and $ (B_t)_{t\geq0} $ denotes a $ d- $dimensional Brownian motion discribing the influence of the surrounding medium. Here we are interested in the scaled equation
%	 Under an appropriate scaling (see e.g. \cite[Chapter 2.2.2]{LRS10}) one obtains the following equation
}
\stepcounter{equation}
dX^\varepsilon_t&=\frac{1}{\varepsilon}V^\varepsilon_tdt,\tag{1.2a}\label{Langeps1}\\ dV^\varepsilon_t&=-\frac{1}{\varepsilon}\nabla\Phi_1(X^\varepsilon_t)dt-\frac{1}{\varepsilon^2}V^\varepsilon_tdt+\frac{1}{\varepsilon}\sqrt{2}dB_t\tag{1.2b}\label{Langeps2},\\
\intertext{cp.~e.g.~\cite[Chapter 2.2.2]{LRS10}. Small $ \varepsilon>0 $ represent the \textit{overdamped regime}. Physically this corresponds to large friction forces and an appropriate time-scaling (see \cite{LRS10}[Chapter 2.2.4] for a physical interpretation). The authors of \cite{RXZ17} prove convergence in law of $ \left(X^\varepsilon_t\right)_{t\geq0} $ as $ \varepsilon $ tends to zero to a solution of the \textit{overdamped Langevin} equation
}
\stepcounter{equation}
dX^0_t&=-\nabla\Phi_1(X^0_t)dt+\sqrt{2}dB_t.\tag{1.3}\label{Overdamped1}%\\
\intertext{Depending on the context a solution to (\ref{Overdamped1}) is also called a distorted Brownian motion. This convergence is known as the \textit{overdamped limit}. More generally, we treat a scaling limit of \textit{generalized stochastic Hamiltonian systems} (gsHs), i.e.,
}
\stepcounter{equation}
\stepcounter{equation}
dX^\varepsilon_t&=\frac{1}{\varepsilon}\nabla\Phi_2(V^\varepsilon_t)dt,\tag{1.4a}\label{sHseps1}\\ dV^\varepsilon_t&=-\frac{1}{\varepsilon}\nabla\Phi_1(X^\varepsilon_t)dt-\frac{1}{\varepsilon^2}\nabla\Phi_2(V^\varepsilon_t)dt+\frac{1}{\varepsilon}\sqrt{2}dB_t.\tag{1.4b}\label{sHseps2}
\end{align}
\noindent Here $ \Phi_2 $ is a potential, generalizing the kinetic energy of the particles, i.e., the Hamiltonian is given by $ \text{H}_{\mathbf \Phi}(x,v)=\Phi_1(x)+\Phi_2(v) $. Observe that for $ \Phi_2(v)=\frac{1}{2}\abs{v}^2 $ we just recover (\ref{Langeps1}),  (\ref{Langeps2}). The main result of this paper is to prove convergence in law of the positions $ \left(X_t^\varepsilon\right)_{t\geq0} $ of (\ref{sHseps1}), (\ref{sHseps2}) to $ \left(X_t^0\right)_{t\geq0} $ from (\ref{Overdamped1}) as $ \varepsilon\to0 $.
%The scaling from (\ref{Langeps1}), (\ref{Langeps2}) can be described in terms of a scaled velocity potential $ \Phi_2^\varepsilon $ in the sense that the corresponding generators are unitarily equivalent. 
 Our assumptions on $ \Phi_1$ and $\Phi_2 $ are so weak that standard results on existence do not apply, see in particular Assumption \ref{C1} and \ref{C2} below. Furthermore, our assumptions allow singular pair interactions like the Lennard-Jones potential. For the pair $ \mathbf \Phi=(\Phi_1,\Phi_2) $ we prove existence of weak solutions $ \left(X_t^\varepsilon,V_t^\varepsilon\right)_{t\geq0} $ to (\ref{sHseps1}), (\ref{sHseps2}) via martingale solutions $ \mathbb P^\varepsilon_{\mathbf \Phi} $ to the generator $ L^\varepsilon_{\mathbf \Phi} $ of (\ref{sHseps1}), (\ref{sHseps2}) given through It\^{o}'s formula, i.e.,
\begin{equation}\label{Lbg}
L^\varepsilon_{\mathbf \Phi}f=\frac{1}{\varepsilon^2}\left(\Delta_vf-\nabla_v\Phi_2\cdot\nabla_vf\right)+\frac{1}{\varepsilon}\left(\nabla_v\Phi_2\cdot\nabla_xf-\nabla_x\Phi_1\cdot\nabla_vf\right) 
\end{equation}
for $ f\in C^\infty_c(\{\Phi_1,\Phi_2<\infty\})$. Observe that the linear operator fails in general to be sectorial, due to the degeneracy of the Laplacian. Hence, the corresponding operator semigroups are not analytic, which males the analysis more challenging.
 
As an intermediate step we consider 
%For $ \varepsilon>0 $ we consider 
for the scaled velocity potential $ \Phi_2^\varepsilon(\cdot)=\Phi_2(\frac{\cdot}{\varepsilon})+\ln{(\varepsilon^d)} $ the pair of potentials $ \mathbf \Phi^\varepsilon=(\Phi_1,\Phi_2^\varepsilon) $. 
%Observe that for $  \Phi_2(v)=\frac{1}{2}\abs{v}^2 $ this leads to a different stochastic differential equation than (\ref{Langeps1}), (\ref{Langeps2}).
The major challenge is to prove weak convergence of the position marginals $ \mathbb P^{1,X}_{\mathbf \Phi^\varepsilon} $ of martingale solutions $ \mathbb P^1_{\mathbf \Phi^\varepsilon} $ corresponding to $ L^1_{\mathbf \Phi^\varepsilon} $ as $ \varepsilon\to 0 $. This we achieve with analytic and probabilistic methods.
The analytic part consists of a semigroup convergence result, the probabilistic one of a tightness result.
At the end we use this convergence and unitary transformations to show convergence of the positions of (\ref{sHseps1}), (\ref{sHseps2}) to a distorted Brownian motion. 
% based on a non-sectorial Lyons-Zheng decomposition. 
%  This result confirms the convergence of the positions $ \left(X^\varepsilon_t\right)_{t\geq0} $ of the scaled Langevin dynamics even for singular interacting forces $ \nabla\Phi_1 $.
% The cost of our tightness result is that we restrict to the invariant measure of the generators as initial distribution of our martingale solutions.

The organization of this paper is as follows.
In Section 2 and 3 we closely follow the approach in \cite{CG10} where martingale solutions for $ \Phi_2=\frac{1}{2}\abs{v}^2 $ were constructed.
Section 2 contains essential m-dissipativity results for the generator $ \left(L^1_{\mathbf \Phi} ,C^\infty_c(\{\Phi_1,\Phi_2<\infty\})\right)$ on $ L^2(\mu_{\mathbf \Phi}) $ and $ L^1(\mu_{\mathbf \Phi}) $, where $ \mu_{ \mathbf\Phi} $ is an invariant measure for $ L^1_{\mathbf{\Phi}} $ from (\ref{Lbg}). In Section 3 we show existence of a martingale solution to $ L^1_{\mathbf{\Phi}} $ in terms of a right process. Section 4 gives a brief overview of the functional analytic objects corresponding to the overdamped Langevin equation (\ref{Overdamped1}) and existence of martingale solutions for its generator is shown. 
The analytic part for convergence is provided in Section 5. We prove strong convergence of the semigroups generated by the scaled generators $ L^1_{\mathbf\Phi^\varepsilon} $. Note that for each $ \varepsilon>0 $ the generator $ L^1_{\mathbf\Phi^\varepsilon} $ is acting on a different Hilbert spaces. Hence, we use the concepts developed by Kuwae--Shioya in \cite{KS03} for showing convergence. Section 6 contains the probabilistic part for convergence. We establish convergence in law of weak solutions via semigroup convergence and tightness of the family $ \left(\mathbb P^1_{\mathbf \Phi^\varepsilon}\right)_{\varepsilon>0} $.
% based on a non-sectorial Lyons-Zheng decomposition.
In Section 7 we explain how these results apply to the original problem, i.~e.~ to prove convergence in law of the positions $ \left(X^\varepsilon_t\right)_{t\geq0} $ from (\ref{sHseps1}), (\ref{sHseps2}) towards $ \left(X^0_t\right)_{t\geq0} $ from (\ref{Overdamped1}).
The core results achieved in this paper may be summarized in the following list:
\begin{itemize}
	\item We prove that the closure of $ \left(L^1_{\mathbf \Phi},C^\infty_c(\{\Phi_1,\Phi_2<\infty\})\right) $ in $ L^1(\mu_{ \mathbf\Phi}) $ is the generator of a sub-Markovian strongly continuous contraction semigroup $ \sg{\mathbf{\Phi}}{,1} $, see Theorem \ref{mdissiL1}.
	\item For the scaled velocity potential $ \Phi_2^\varepsilon $ we prove convergence of the associated $ L^2(\mu_{\mathbf{\Phi}^\varepsilon}) $ semigroups $ \sg{\mathbf{\Phi}^\varepsilon}{,2} $ in the sense of Kuwae--Shioya, see Theorem \ref{CoSe}.
	\item We prove weak convergence of the position marginals $\mathbb P^{1,X}_{\mathbf \Phi^\varepsilon} $, ,$ \varepsilon>0 $, to a martingale solution of the generator of the distorted Brownian motion as $ \varepsilon\to0 $, 
%	for the generators $ L_{\mathbf{\Phi}^\varepsilon} $ via a non-sectorial Lyons-Zheng decomposition,
	 see Corollary \ref{FinalConResult}.
	\item We give a rigorous proof for the convergence in law of the positions $ \left(X_t^\varepsilon\right)_{t\geq0} $ of weak solutions $ \left(X_t^\varepsilon,V_t^\varepsilon\right)_{t\geq0} $ to (\ref{sHseps1}), (\ref{sHseps2}) to the overdamped Langevin equation as $ \varepsilon\to0 $, see Theorem \ref{finalresult}.
\end{itemize}
At this point we would like to point out that all results hold for very large class of interaction potentials $ \Phi_1 $ which can also be very singular, e.g., potentials of Lennard--Jones type are admissible.

Our results are complementary to those in \cite{RXZ17} in the following sense:
First, there the authors have to assume the interaction term $ \nabla\Phi_1 $ to be continuous. Second, there the state space is assumed to be the $ d- $dimensional torus $ \mathbb T^d $. Due to our weaker assumptions the weak solutions constructed in our framework require initial distributions which are absolutely continuous w.r.t.~ the invariant measure $ \mu_{ \mathbf\Phi} $. This aspect is more restrictive than in \cite{RXZ17}. 
Additionally, the $ \Phi_1 $ in \cite{RXZ17} may also depend on $ \varepsilon>0 $. 
%Assuming global Lipschitz continuity of $ \nabla\Phi_1 $, pointwise in time $ L^2-$convergence of strong solutions is shown in \cite{SMD82}.
% !TeX spellcheck = en_US

\section[M-Dissipativity of the Operator $L^1_{\mathbf{\Phi}}$]{M-Dissipativity of the Operator $L^1_{\mathbf{\Phi}}$}\label{Mdiss}
The main goal of this section is to establish for a pair $ \mathbf{\Phi}=(\Phi_1,\Phi_2) $ of potentials essential m-dissipativity of the differential operator
$ \left(L^1_{\mathbf{\Phi}},C^\infty_c(\{\Phi_1,\Phi_2 <\infty\})\right) $ given by
\begin{equation}\label{L00}
L^1_{\mathbf{\Phi}}f=\Delta_vf-\nabla_v\Phi_2\cdot\nabla_vf+\nabla_v\Phi_2
\cdot\nabla_xf-\nabla_x\Phi_1\cdot\nabla_vf, \quad f\in C^\infty_c(\{\Phi_1,\Phi_2 <\infty\})
\end{equation}
on $ L^1(\Rd{2d},\mu_\mathbf{\Phi}) $, where $\mu_\mathbf{\Phi}$ is
% an invariant measure for $ L_{\mathbf{\Phi}} $ and 
 absolutely continuous w.r.t.~ the Lebesque measure on $ \left(\Rd{2d},\mathcal B(\Rd{2d})\right) $. In the following we always denote $ L^1_{\mathbf{\Phi}} $ by $ L_{\mathbf{\Phi}} $.
%  The constants $ \beta,\gamma $ from (\ref{Lbg}) can be assumed to be 1 (see \cite{C10}[p. 151]).
  We follow closely the argumentation in \cite{CG10} and generalize the proofs therein for a general velocity potential $ \Phi_2 $
 fulfilling the Assumptions \ref{C2} below. Therefore we only prove the parts which actually differ and refer to \cite{CG10} for additional details. First we prove essential m-dissipativity on $L^2(\Rd{2d},\mu_\mathbf{\Phi}) $ for locally Lipschitz continuous $ \Phi_1 $.
Afterwards we use this result to show the m-dissipativity of the closure of (\ref{L00}) on $ L^1(\Rd{2d},\mu_\mathbf{\Phi})  $ for singular $ \Phi_1 $.
The potentials $ \Phi_1 $, $ \Phi_2 $ and their derivatives are considered as functions on $ \Rd{2d} $ and $ \Rd{d} $
simultaneously in the following way:
$\Phi_1(x,v)=\Phi_1(x)$, $ \Phi_2(x,v)=\Phi_2(v) $, where $ (x,v)\in\Rd{d}\times\Rd{d} $. For a (weakly) differentiable function $ f $ on $ \Rd{2d}
 $, $ \nabla_x f $ denotes the $ d- $dimensional (weak) gradient w.r.t. the first $ d $ unit vectors. 
 Corresponding definitions hold for $\nabla_v,\Delta_x,\Delta_v,\partial_{x_i},\partial_{v_i} $, $ i=1,...,d$. Expression like $ \nabla_v\Phi_2\cdot\nabla_vf $ from (\ref{L00}) are understood as $  \nabla_v\Phi_2\cdot\nabla_vf(x,v)=\sum_{i=1}^d \partial_{v_i}\Phi_2(x,v)\partial_{v_i}f(x,v)$.
 The gradient, the Laplacian and weak partial derivatives of $ \Phi_1 $ and $ \Phi_2 $ considered as a function on $\Rd{d}$ are denoted by $  \nabla , \Delta ,\partial_i$,  $i=1,...,d $, respectively.
% \begin{D}
% 	% \mbox{}\vspace{-\topskip}
% 	Let $\left(L,D(L)\right)$ be a linear operator on a Hilbert space $\mathcal{H}$.
% 	$\left(L,D(L)\right)$ is called 
% 	\begin{enumerate}
% 		\item[(i)] \textit{dissipative} if 
% 		\begin{equation*}
% 		\Re(Lf,f)_\mathcal{H}\leq 0, \quad\forall f\in D(L).
% 		\end{equation*}
% 		\item[(ii)] \textit{m-dissipative} if $\left(L,D(L)\right)$ is densely defined, dissipative and $\mathcal{R}(\alpha-L)=\mathcal{H}$ for one $\alpha>0$.
% 		\item[(iii)] \textit{essentially m-dissipative} if $\left(L,D(L)\right)$ is densely defined, dissipative
% 		and $\mathcal{R}(\alpha-L)$ is dense in $\mathcal{H}$ for one $\alpha>0$.
% 	\end{enumerate}
% \end{D}
%The importance of m-dissipative operators is due to the Lumer-Phillips theorem \cite[Theorem 3.1]{LP61}:
% \begin{T}[Lumer-Philips.]
%	Suppose $(L, D(L))$ is a linear operator on $\mathcal{H}$. Then $(L,D(L))$ is the
%	generator of a strongly continuous contraction semigroup on $\mathcal{H}$ if and only if $(L,D(L))$ is m-dissipative.
%\end{T}
\begin{N}\label{mu}
	For $ n\in\mathbb N $ and a measurable function $\Psi:\Rd{n}\longrightarrow \overline{\Rd{}}$, where $ \overline{\Rd{}} $ denotes the extended real numbers, we define the
	measure $\mu_\Psi$ by its Radon-Nikodym derivative w.r.t. the Lebesgue measure $ dx $ on 
	$(\Rd{n},\mathcal{B}(\Rd{n}))$, i.e.,
	\begin{equation*}
	\frac{d\mu_\Psi}{dx}=e^{-\Psi}.
	\end{equation*}
\end{N}
\noindent We state the assumptions we later assume for the position potential $ \Phi_1 $ and the velocity potential $ \Phi_2 $: 
\begin{C}\label{C1}
	Let $\Phi_1:\Rd{d}\longrightarrow\mathbb{R}\cup \{\infty\}$ and $ q\in[2,\infty] $.
	\begin{itemize}[align=parleft, labelsep=0.5cm]
		\item [$(\Phi_1\thePHIeins)$]
		\stepcounter{PHIeins}
		\tabto{1cm}
		$\Phi_1$ is locally Lipschitz continuous, i.e.,
		the restriction of $\Phi_1$ to an arbitrary \tabto{1cm}
		compact subset of $\Rd{d}$ is Lipschitz continuous. In particular, $\Phi_1:\Rd{d}\longrightarrow\mathbb{R}$.
		\item [$(\Phi_1\thePHIeins)$]
		\stepcounter{PHIeins}
		\tabto{1cm}
		$ \Phi_1 $ is bounded from below and $\{\Phi_1<\infty\}\neq \emptyset $.
		\item [$(\Phi_1\thePHIeins)$]
		\stepcounter{PHIeins}
		\tabto{1cm} $e^{-\Phi_1} $ is continuous on $\Rd{d} $.
		\item [$(\Phi_1\thePHIeins)^q$]
		\stepcounter{PHIeins}
		\tabto{1cm} $ \Phi_1 $ is weakly differentiable on $ \{\Phi_1<\infty\} $ and 
		$ \nabla \Phi_1\in L^q_{loc}(\Rd{d},\mu_{\Phi_1}) $.
	\end{itemize}
\end{C}

\begin{C}\label{C2}
	Let $\Phi_2:\Rd{d}\longrightarrow\mathbb{R}\cup \{\infty\}$.
	\begin{itemize}[align=parleft, labelsep=0.5cm]
		\item [$(\Phi_2\thePHIzwei)$]
		\tabto{1cm}
		$\Phi_2$ is $\mathcal{B}(\Rd{d})-\mathcal{B}(\overline{\mathbb{R}})$ measurable and $\{\Phi_2 <\infty\}\neq\emptyset$ is open.
		\stepcounter{PHIzwei}		     
		\item [$(\Phi_2\thePHIzwei)$]\label{Phi2}\tabto{1cm}
		 $\Phi_2$ is bounded from below and locally integrable on $ \{\Phi_2<\infty\} $.
		\stepcounter{PHIzwei}
		\item [$(\Phi_2\thePHIzwei)$]\tabto{1cm}
		For $i\in\{1,..,d\}$ it holds for the distributional derivatives\tabto{1cm}
		{$ \partial_i\Phi_2\in L_{loc}^2(\{\Phi_2 <\infty\})$} and $
		\partial_i^2\Phi_2\in L_{loc}^1(\{\Phi_2 <\infty\}) $.
		%			$\partial_i^2\Phi_2\in L_{loc}^2(\{\Phi_2 <\infty\})$, 
		%		$\partial_i^2\Phi_2\in\mathcal{D}^\prime(\{\Phi_2 <\infty\})
		%		\cap L_{loc}^2(\{\Phi_2 <\infty\})$. 
		\stepcounter{PHIzwei}
		\item [$(\Phi_2\thePHIzwei)$]\tabto{1cm}
		$\left(\Delta-\nabla\Phi_2\cdot\nabla,
		C^\infty_c(\{\Phi_2 <\infty\})\right)$ is essentially self-adjoint 
		on $L^2(\Rd{d},\mu_{\Phi_2})$.	    
		\stepcounter{PHIzwei}
		\item[$(\Phi_2\thePHIzwei)$]\tabto{1cm} There are constants $K\in (0,\infty)$ and $\alpha\in[1,2)$ such that it holds\tabto{1cm} {$\abs{\Delta\Phi_2}\leq
			K(1+\abs{\nabla\Phi_2}^\alpha)$.}
		\stepcounter{PHIzwei}
	\end{itemize}
\end{C}
\noindent According to Notation \ref{mu} denote by $ \mu_{\mathbf\Phi}$ the measure $\mu_{\Phi_1+\Phi_2} $ on $ \left(\Rd{2d},\mathcal B(\Rd{2d})\right) $ and by $ \mathcal{H}_{\mathbf\Phi}$ the Hilbert space $L^2(\Rd{2d},\mu_{\mathbf\Phi}) $. 
%In the following we assume that $\Phi_1$ and $\Phi_2$ 
%fulfill Assumption \ref{C1} and \ref{C2}, respectively. 
\begin{Rem}\label{densedefined}
	\mbox{}\vspace{-\topskip}	
	~\begin{enumerate}
		\item[(i)]\label{LipHinf}
		Let $\Omega$ be an open subset of $\Rd{d}$. Then it holds $f\in H^{1,\infty}_{loc}(\Omega)$ if and only if $f$ has a representative which
		is locally Lipschitz continuous in $\Omega$ (see \cite[Chapter 5.8, Theorem 4]{Ev10}). Hence, the assumption $(\Phi_11) $ implies $(\Phi_12)-(\Phi_14)^\infty$
		apart from the boundedness from below.
		\item[(ii)] If we assume instead of $ (\Phi_22) $ the following condition:
		\begin{itemize}[align=parleft, labelsep=0.5cm]
			\item[$\widetilde{(\Phi_22)}$]$\Phi_2 $ is locally bounded on $ \{\Phi_2 <\infty\} $.
		\end{itemize}
		Then in combination with $ (\Phi_25) $ one can argue similar as in the proof of \cite{C10}[Lemma A6.2.] that $ \Phi_2 $ is continuously differentiable on $ \{\Phi_2 <\infty\} $ and $ \nabla\Phi_2 $ is locally Lipschitz on $ \{\Phi_2 <\infty\} $.
		\item[(iii)]	Assuming $ (\Phi_12),(\Phi_14)^q,(\Phi_22)$ and  $(\Phi_23) $ we can consider $ \left(L_{\mathbf{\Phi}},C^\infty_c(\{\Phi_1,\Phi_2 <\infty\})\right) $ as an operator on $ L^p(\Rd{2d},\mu_{ \mathbf\Phi}) $ for every $ p\in [1,2] $.
		\item[(iv)]\label{densesubspace}
		Since the measure $ \mu_{\Phi_2} $ on $ \Rd{d} $ is locally finite it holds by \cite[Proposition 7.2.3]{Co80} that $ \mu_{\Phi_2} $ is regular Borel measure on $ \left (\{\Phi_2 <\infty\},\mathcal B(\{\Phi_2 <\infty\}) \right )$ and hence by \cite[Proposition 7.4.2]{Co80} the set
		$C^\infty_c(\{\Phi_2 <\infty\})$ is dense in $L^2(\{\Phi_2<\infty\},\mu_{\Phi_2})\cong L^2(\Rd{d},\mu_{\Phi_2})$. 		
		%		This can be seen by approximating an arbitrary element from $ L^2(\Rd{d},\mu_{\Phi_2}) $ first by elements which vanish outside some compact set and employing \cite[Corollary 2.30]{AF08} afterwards.
		%%		.\newline 
		%		Since $\mu_{\Phi_2}(\{\Phi_2=\infty\})=0$ it is enough to prove 
		%		that $C^\infty_c(\{\Phi_2 <\infty\})$ is dense in $L^2(\{\Phi_2<\infty\},\mu_{\Phi_2})$.
		%		Since $\Phi_2$ is assumed to be locally Lipschitz continuous in $\{\Phi_2<\infty\}$,
		%		$\mu_{\Phi_2}$ assigns finite measure to compact subsets of $\{\Phi_2 <\infty\}$.
		%		\cite[Proposition 7.2.3]{Coh13} implies therefore that $\mu_{\Phi_2}$ is regular.
		%		Therefore we can apply
		%		\cite[Proposition 7.4.2]{Coh13} and obtain 
		%		that the space of continuous functions with compact support $C_c(\{\Phi_2 <\infty\})$ 
		%		is dense in $L^2(\{\Phi_2<\infty\},\mu_{\Phi_2})$.
		%		Thus, it suffices to show that continuous and compactly supported functions can be approximated by elements from  
		%		$C^\infty_c(\{\Phi_2 <\infty\})$. So let $0\neq f\in C_c(\{\Phi_2 <\infty\})$.
		%		We choose a cutoff function $\eta\in C^\infty_c(\{\Phi_2 <\infty\})$ on supp$(f)$ such that $0\leq\eta\leq1$.
		%		By the Weierstra{\ss} approximation theorem we can find a sequence of polynomials $p_n$ converging uniformly 
		%		on supp$(\eta)$ to $f$. The functions $f_n=\eta p_n$ are obviously elements from 
		%		$C^\infty_c(\{\Phi_2 <\infty\})$ and one easily checks $f_n\longrightarrow f$ in 
		%		$L^2(\{\Phi_2<\infty\},\mu_{\Phi_2})$ as $n\to\infty$ because $\mu_{\Phi_2}$ is locally
		%		finite.
		\item[(v)] See Remark \ref{RemSingPot} as a reference for sufficient conditions implying $(\Phi_24)$.
	\end{enumerate}
\end{Rem}

\begin{Prop}\label{H1infty}
	Let $\Omega\subseteq\Rd{n}$, $ n\in\mathbb N $, be open and $\Psi:\Omega\longrightarrow \Rd{}$ be measurable and locally bounded or bounded from below and locally integrable. Assume further that the first order distributional derivatives $ \partial_i\Psi $, $ i\in\{1,...,n\} $, are in 
	$ L^p_{loc}(\Omega) $, for some $ p\in[1,\infty] $. Then it holds that $e^{-\Psi}\in H^{1,p}_{loc}(\Omega)$ and $\partial_i\left(e^{-\Psi}\right)=-\partial_i\Psi e^{-\Psi}$.
	\begin{proof}
		Let $ \Omega^\prime\subset\Omega $ be open such that $ \overline{\Omega^\prime}\subseteq\Omega $ is compact. We need to show that $e^{-\Psi}\in H^{1,p}(\Omega^\prime)$. Hence, let
		$\varphi\in C^\infty_c(\Omega^\prime)$ be arbitrary.	
		Since $K:=\text{supp}(\varphi)$ is compact there is a non-negative $ \chi \in C^\infty_c(\Omega^\prime)$ such that $ \chi=1 $ on $ K $. Obviously $e^{-\Psi}\in L^\infty(\Omega^\prime)\subseteq L^p(\Omega^\prime)$.
		By the compact support of $ \chi $ and a regularization as in \cite[Lemma 3.16]{AF08} one can find a sequence $ \left(u_k\right)_{k\in\mathbb{N}}\in C^\infty_c(\Omega^\prime) $ such that  $u_k\longrightarrow \chi\Psi$, as $ k\to\infty, $ in $H^{1,1}(\Omega^\prime)$. 
		In the case of locally bounded $ \Psi $ it holds $\norm{u_k}_\infty\leq\norm{\chi\Psi}_\infty$, for all $ k\in\mathbb N$. Otherwise, if $ C\in \Rd{} $ is a lower bound of $ \Psi $ then it holds
		$ C\leq u_k(x) $ for all $ x\in\Omega^\prime $ and all $ k\in\mathbb N $.  
		By switching to a subsequence which we also denote by $ (u_k)_{k\in\mathbb N} $ we can apply the dominated convergence theorem, integration by parts and H\"olders inequality to obtain
		\begin{align*}
		\int\displaylimits_{\Omega^\prime} e^{-\Psi}\partial_i\varphi\,dx=\lim_{k\to\infty}\int\displaylimits_{\Omega^\prime} e^{-u_k}\partial_i\varphi\,dx
		=\lim_{k\to\infty}\int\displaylimits_{\Omega^\prime} \partial_iu_ke^{-u_k}\varphi\,dx
		=\int\displaylimits_{\Omega^\prime}\partial_i\Psi e^{-\Psi}\varphi\,dx.
		\end{align*}
	\end{proof}
\end{Prop}
\noindent Under the assumptions $(\Phi_12)-(\Phi_14)^q $, $ q\in[2,\infty] $ and $(\Phi_21)-(\Phi_23)$ we obtain the following proposition and corollary:
\begin{Prop}\label{SymAntisy}
	$\left(L_{\mathbf\Phi},C_c^{\infty}(\{\Phi_1,\Phi_2<\infty\})\right)$ admits a decomposition into 
	$ L_{\mathbf{\Phi}}=S+A $, with  symmetric $ S $ and antisymmetric $ A $ on $ C_c^{\infty}(\{\Phi_2<\infty\})$ w.r.t. the scalar product on $\mathcal H_{\mathbf\Phi} $. $ S $ and $ A $ are given through 
	\begin{equation*}
	Sf=\Delta_vf-\nabla_v\Phi_2\cdot\nabla_vf,\quad 
	Af=\nabla_v\Phi_2\cdot\nabla_xf-\nabla_x\Phi_1\cdot\nabla_vf,\quad  f\in C_c^{\infty}(\{\Phi_1,\Phi_2<\infty\}).
	\end{equation*}
	\begin{proof}
		The proof consists of the product rule for Sobolev functions and Proposition \ref{H1infty}.			
	\end{proof}
\end{Prop}
\begin{Co}\label{invariant}
	The measure $\mu_{\mathbf\Phi}$ is invariant for 
	$\left(L_{\mathbf{\Phi}},C_c^{\infty}(\{\Phi_1,\Phi_2<\infty\})\right)$, i.e., $  L_{\mathbf{\Phi}}f $ is integrable w.r.t. $ \mu_{\mathbf\Phi} $ for all ${f\in C_c^{\infty}(\{\Phi_1,\Phi_2<\infty\})}$
	and it holds
	\begin{equation}\label{invariant1}
	\int\displaylimits_{\mathbb{R}^{2d}}L_{\mathbf{\Phi}}f\,d\mu_{\mathbf\Phi}=0.
	\end{equation}
%	In particular, $\left(L_{\mathbf{\Phi}},C_c^{\infty}(\{\Phi_1,\Phi_2<\infty\})\right)$ is  dissipative .
	In particular, $\left(L_{\mathbf{\Phi}},C_c^{\infty}(\{\Phi_1,\Phi_2<\infty\})\right)$ is closable and its closure $\left(L_{\mathbf{\Phi},p},D(L_{\mathbf{\Phi},p})\right)$ is dissipative on $L^p(\Rd{2d},\mu_{\mathbf{\Phi}})$ for every $ p\in [1,2] $.
	%	 for all ${f\in C_c^{\infty}(\{\Phi_2<\infty\})$ it holds $ L_{\mathbf{\Phi}}f\in L^1(\Rd{2d},\mu_{\mathbf\Phi})  $ and it holds
	%	\begin{equation}
	%	\int\displaylimits_{\mathbb{R}^{2d}}L_{\mathbf{\Phi}}f\,d\mu_{\mathbf\Phi}=0.
	%	\end{equation}
	\begin{proof} For $ f\in C^\infty_c(\{\Phi_1,\Phi_2<\infty\}) $ one chooses a cut off function $\eta\in C^\infty_c(\{\Phi_1,\Phi_2<\infty\}) $, s.t. $ \eta=1 $ on $ \text{supp}(f) $ and uses the decomposition from Proposition \ref{SymAntisy}. But $ S\eta,A\eta $ vanish on $ \text{supp}(f) $ which implies (\ref{invariant1}). The dissipativity follows by \cite[Lemma 1.8, App. B]{Eb99}. 
	\end{proof}
\end{Co}

\subsection{M-Dissipativity for locally Lipschitz continuous $ \Phi_1 $ on $ L^2(\Rd{2d},\mu_{\mathbf\Phi}) $}
Throughout this first part we assume that $ \Phi_1 $ and $\Phi_2 $ fulfill $ (\Phi_11) $ and $ (\Phi_21)-(\Phi_25) $, respectively. In particular, it holds $ \{\Phi_1<\infty\}=\Rd{d} $.
%\begin{Co}\label{dissi}
%	$\left(L_{\mathbf{\Phi}},C_c^{\infty}(\{\Phi_2<\infty\})\right)$ is a dissipative
%	operator on $\mathcal H_{\mathbf{\Phi}}$.
%	\begin{proof}
%		In the previous corollary we proved that $ \mu_{\mathbf\Phi} $ is an invariant measure. One easily checks that $\left(L_{\mathbf{\Phi}},C_c^{\infty}(\{\Phi_2<\infty\})\right)$ fulfils \cite[Definition 1.5, App. B]{Eb99}. Then \cite[Lemma 1.8, App. B]{Eb99} finishes the proof.
%	\end{proof}	
%\end{Co}
%The following lemma is well-known.
%\begin{Le}\label{Noext}
%	An m-dissipative operator $ (L,D(L)) $ on a Hilbert space $ H $ has no proper dissipative extension.
%\end{Le}
\begin{Prop}\label{sm}
	Let $\left(L,D\right)$ be a densely defined operator 
	on a Hilbert space $\mathcal H$. 
	Furthermore $L$ is assumed to be symmetric and negative definite.
	If $\left(L,D\right)$ is essentially self-adjoint, then $\left(L,D\right)$ is
	essentially m-dissipative.
	\begin{proof}
		Since $ (L,D) $ is negative definite its closure $\left(\bar{L},D(\bar{L})\right)$ is dissipative, implying that $ 1-\bar L $ is injective. By assumption it holds $ \mathcal{R}(1-\bar L)^\perp=\mathcal N(1-\bar L) =\{0\}$.		
	\end{proof}
\end{Prop}
\begin{T}\label{maintheorem}
	Assume $ (\Phi_11) $ and $ (\Phi_21)-(\Phi_25) $. Then the operator $\left(L_{\mathbf{\Phi}},C_c^{\infty}(\{\Phi_2<\infty\})\right)$ is essentially 
	m-dissipative on $\mathcal H_{\mathbf{\Phi}}$. The strongly continuous contraction semigroup $ \left(T^{\mathbf\Phi}_t\right)_{t\geq 0}$ generated by the closure of
	$\left(L_{\mathbf{\Phi}},C_c^{\infty}(\{\Phi_2<\infty\})\right)$
%	 on $\mathcal H_{\mathbf{\Phi}}$
	  is sub-Markovian.
	% $\left(T^{\Phi_1,\Phi_2}_t\right)_{t\geq 0}$.
	\begin{proof}
		This proof is based on the idea of the proof of \cite[Thm. 2.1]{CG10}. In the first part $\Phi_1$ is considered to be globally Lipschitz continuous with Lipschitz constant $C_{\Phi_1}$. The second part treats the general case.
		Throughout the first part of the proof all function spaces consist of complex valued functions. Observe that those spaces are isometric to the complexification of the real valued function spaces. Furthermore, $ L_{\mathbf{\Phi}} $ leaves the real valued functions invariant. Hence, we show that the complexified operator
		is essentially m-dissipative, this proves the theorem for the real cases.
		
		% Lets start with the first part.\\
	\noindent\textbf{1st part:}\\
		The basic idea is to use the unitary transformation 
%		and proof
%		essential m-dissipativity of the transformed operator. We will use the unitary transformation
		\begin{equation}\label{Uni}
		U:L^2(\mathbb{R}^{2d} ,\mu_{\mathbf\Phi})\longrightarrow 
		L^2(\{\Phi_2<\infty\}),\quad
		f\mapsto\exp{(-\frac{\Phi_1+\Phi_2}{2})}f. 
		\end{equation}
%		Assuming $\Phi_1$ and $\Phi_2$ to be $C^\infty$, one can easily check by calculating derivatives
		Formally $\left(L_{\mathbf{\Phi}},C_c^{\infty}(\{\Phi_2<\infty\})\right)$ transforms under $ U $ into the operator 
		\begin{equation}\label{aim}
		L=U L_{\mathbf{\Phi}}U^*=\Delta_v+\frac{\Delta_v\Phi_2}{2}-
		\frac{\abs{\nabla_v\Phi_2}^2}{4}
		+\nabla_v\Phi_2\cdot\nabla_x-\nabla_x\Phi_1\cdot\nabla_v.
		\end{equation}
%		defined on $C_c^{\infty}(\{\Phi_2<\infty\})$.
		In the following we prove essential m-dissipativity of $L$ on a 
		suitable chosen domain $D$. Afterwards we make the transformation in (\ref{aim}) rigorous.
%		 to obtain the desired result for 	$\left(L_{\mathbf{\Phi}},C_c^{\infty}(\{\Phi_2<\infty\})\right)$. 
		 Assumption $(\Phi_24)$ gives us the negative definite and essentially self-adjoint operator
		$\left(\Delta-\nabla\Phi_2\cdot\nabla,C_c^\infty(\{\Phi_2<\infty\})\right)$
		 on $L^2(\mathbb{R}^{d} ,\mu_{\Phi_2})$. 
		Proposition \ref{sm} implies that $\left(\Delta-\nabla\Phi_2\cdot\nabla,C_c^{\infty}(\{\Phi_2<\infty\})\right)$ is essentially m-dissipative on $L^2(\Rd{d},\mu_{\Phi_2})$. Consider the unitary transformation
		\begin{equation}\label{unitary1}
		U_{\Phi_2}:L^2(\mathbb{R}^{d} ,\mu_{\Phi_2})\longrightarrow
		L^2(\{\Phi_2<\infty\}),\quad g\mapsto\exp{(-\frac{1}{2}\Phi_2)}g .
		\end{equation}
		Since unitary transformations preserve essential m-dissipativity we have that
		\begin{equation}\label{L1}
		L_0=U_{\Phi_2} (\Delta-\nabla\Phi_2\cdot\nabla)U^*_{\Phi_2}
		\end{equation}
		defined on $U_{\Phi_2} C_c^{\infty}(\{\Phi_2<\infty\})$  
		is an essentially m-dissipative operator on $L^2(\{\Phi_2<\infty\})$. Let $ g\in C_c^{\infty}(\{\Phi_2<\infty\}) $ and $ f=U_{\Phi_2}g $. In the following the differential operators $ \Delta$ and $\nabla $ are understood in the distributional sense. Then it holds
		\begin{equation}\label{distlap}
		\Delta f=\Delta (U_{\Phi_2}g)=\Delta g \exp(-\frac{1}{2}\Phi_2)+2\nabla \left(\exp(-\frac{1}{2}\Phi_2)\right )\cdot\nabla g+g\Delta \exp(-\frac{1}{2}\Phi_2).
		\end{equation}
		Proposition \ref{H1infty} and (\ref{distlap}) lead to
%		Due to the product rule for distributions and Proposition \ref{H1infty} combined with the assumptions in $({\Phi}_23) $ we obtain for $L_0$ the following representation:
		\begin{align}
			L^2(\{\Phi_2<\infty\})\ni L_0 f&=U_{\Phi_2} (\Delta-\nabla\Phi_2\cdot\nabla)g\label{L02}\\
			&=\Delta g\exp(-\frac{1}{2}\Phi_2)+2\nabla \left(\exp(-\frac{1}{2}\Phi_2)\right )\cdot\nabla g\notag\\
			&=\Delta f-g\Delta \exp(-\frac{1}{2}\Phi_2)\label{L01}			
		\end{align} 
	Due to  the Assumptions in $({\Phi}_23) $ and an approximation procedure as in the proof of Proposition \ref{H1infty} one has	
	$\Delta\exp(-\frac{1}{2}\Phi_2)=-\left(\frac{\Delta\Phi_2}{2}-\frac{\abs{\nabla\Phi_2}^2}{4}\right)\exp(-\frac{1}{2}\Phi_2)$, which gives in (\ref{L01})
	\begin{equation}\label{L0}
	 L_0 f=\Delta f+\left(\frac{\Delta\Phi_2}{2}-\frac{\abs{\nabla\Phi_2}^2}{4}\right)f,\text{ for all } f\in U_{\Phi_2}C_c^{\infty}(\{\Phi_2<\infty\}).
	\end{equation}
	Note: The single summands $ \abs{\nabla\Phi_2}^2f $ and $ \Delta\Phi_2f $ in (\ref{L0}) are not necessarily in $ L^2(\{\Phi_2<\infty\}) $. Anyways, $ L_0f $ is an element of $ L^2(\{\Phi_2<\infty\}) $ which can be seen by (\ref{L02}). 
	Nevertheless, $ (\ref{L0}) $ is a suitable representation of $ L_0f $. 	
		Furthermore, $L_0$ is still symmetric and negative definite because we obtained $L_0$ from
		a unitary transformation of a symmetric and negative definite operator.
		
		So far we only worked on the velocity component. To take the position variable $x$ into %account we consider $L_0$ as an operator on $L^2(\{\Phi_2<\infty\})$, i.e., we 
		into account we define a new domain $ D_0\subseteq L^2(\{\Phi_2<\infty\},) $
		\begin{align}
		D_0&:=L^2_c(\mathbb{R}^d)\otimes U_{\Phi_2} C_c^{\infty}(\{\Phi_2<\infty\})\notag\\
		&:=\text{span}\left \{\Rd{2d}\ni(x,v)\mapsto f(x)g(v)\mid f\in L^2_c(\mathbb{R}^d),g\in U_{\Phi_2} C_c^{\infty}(\{\Phi_2<\infty\}) \right \}\label{DefTensor}
		\end{align}
		where $L_c^2(\Rd{d})$ denotes the subspace of $L^2(\Rd{d})$ with elements vanishing	almost everywhere outside a bounded set. For $f=h\otimes g\in D_0$ we set	$L^\prime_0f:=h\otimes  L _0g=\Delta_vf -\frac{\abs{\nabla_v\Phi_2}^2}{4}f+\frac{\Delta_v\Phi_2}{2}f.$
		We extend $L_0^\prime$ linearly to $D_0$. In the following we denote
		the norm and inner product of $L^2(\{\Phi_2<\infty\})$ by $\norm{\cdot}$
		and $(\cdot,\cdot)$, respectively. Let's make some observations on $(L_0^\prime,D_0)$:
		\begin{enumerate}\label{dense}
%			\item[(i)] $(L^\prime_0, D_0)$ is well-defined, i.e. for 
%			$f=\sum_{i=1}^mh_i\otimes g_i\in D_0$, $f=0\Longrightarrow L_0^\prime f=0$.			
%			This holds since multiplication by a function and taking derivatives is a linear
%			operation.
			\item[(i)] $(L_0^\prime,D_0)$ is symmetric, negative definite and densely defined.			
			\item[(ii)] $(L^\prime_0, D_0)$ is essentially m-dissipative.
		\end{enumerate}
		We perturb $L^\prime_0$ with the multiplication operator $\left(B_0,D_0\right)$ given by
%		\begin{equation}
%		B_0f=i\nabla_v\Phi_2\cdot xf,\quad f\in D_0, 
%		\end{equation}
		the measurable function
		\[i\nabla_v\Phi_2\cdot x:\{\Phi_2<\infty\}\longrightarrow\mathbb C,\quad(x,v)\mapsto 
		i\nabla_v\Phi_2(x,v)\cdot x:=i\sum_{l=1}^d\partial_l{\Phi}_2(v)x_l.\]
		Since $\nabla_v\Phi_2\cdot x$ is real valued it follows that $B_0$ is antisymmetric,
%		, i.e. 
%		\[(i\nabla_v\Phi_2\cdot x f,g)=-\overline{(i\nabla_v\Phi_2\cdot x g,f)}\]
%		for $f,g\in D_0$ and
		in particular, $\left(B_0,D_0\right)$ is dissipative.
		We consider the complete orthogonal family of projections $ \left(P_k\right)_{k\in\mathbb N} $ given by 
		\[P_k:L^2(\{\Phi_2<\infty\})\longrightarrow L^2(\{\Phi_2<\infty\}),
		f\mapsto g_kf,\]
		where $g_k(x,v)=1_{[k-1,k]}(\abs{x}_2)$, $ k\in\mathbb{N}$.
		Obviously each $P_k$ maps $D_0$ into itself and $L^\prime_0$ as well as $B_0$ commute with
		each $P_k$ on $D_0 $. In order to apply \cite[Lemma 3]{CG08}
		we need to show that $B_0^k:=P_kB_0$ is $L_k:=P_kL_0^\prime$ bounded with $L_k$-bound less then one. By the Cauchy-Schwarz inequality and the definition of $ P_k $ we have
		\begin{equation}
		\abs{\nabla_v\Phi_2\cdot x}^2\abs{f}^2\leq k^2\abs{\nabla_v\Phi_2}^2\abs{f}^2, \text{ for }f\in P_k D_0.
		\end{equation}
%		 for $f\in P_k D_0$ the following estimate $\abs{\nabla_v\Phi_2\cdot x}^2\abs{f}^2\leq k^2\abs{\nabla_v\Phi_2}^2\abs{f}^2.$
		Hence, it suffices to show that $\norm{\abs{\nabla_v\Phi_2}f}^2\leq a(L^\prime_0f,f)+b\norm{f}^2$ 
		holds for some finite constants $a,b$ independent of $f\in P_k D_0$. Therefore, let $f\in D_0$ and observe that $-\Delta_v$ is positive definite on $D_0$ and $\Delta_v\Phi_2f\in
		L^2(\{\Phi_2<\infty\})$ due to assumption $(\Phi_23)$. Due to the assumptions on $ f $ and $ \Phi_2 $ it holds
		\begin{equation}\label{nablaPhi2}
		\norm{\abs{\nabla_v\Phi_2} f}^2 \leq 4\left(-\Big(\Delta_v-\frac{\abs{\nabla_v\Phi_2}^2}{4}
		+\frac{\Delta_v\Phi_2}{2}\Big)f,f\right)+2(\Delta_v\Phi_2 f,f) 
		\end{equation}
		with both summands on the right-hand side being finite.
		Let $K>0$ and $1\leq\alpha<2$ be the constants from assumption $(\Phi_25)$. Then we have the following estimate
		for the last term in (\ref{nablaPhi2}) 
		\begin{equation}\label{nablaPhizweizwei}
		(\Delta_v\Phi_2 f,f)\leq K\left(\norm{f}^2 
		+\int\displaylimits_{\{\Phi_2<\infty\}}\abs{\nabla_v\Phi_2}^\alpha \abs{f}^{2}\,d(x,v)\right)
		\end{equation}
%		Let $\gamma$ be an arbitrary positive constant.
%		Set $ \gamma=(2\alpha K)^{\frac{\alpha}{2}} $ and apply
		 H\"older's and Young's inequality imply for the last integral on the right hand side
		of $(\ref{nablaPhizweizwei})$ for $p=\frac{2}{\alpha}$, $q=\frac{2}{2-\alpha}$
%		\begin{align}
%		\int\displaylimits_{\{\Phi_2<\infty\}} \abs{\nabla \Phi_2}^\alpha \abs{f}^{2}\,d(x,v)
%%		&=\int\displaylimits_{\{\Phi_2<\infty\}} \abs{\nabla \Phi_2f}^\alpha
%%		\abs{f}^{2-\alpha}\,d(x,v)\notag\\
%%		&\leq \frac{1}{\gamma}\left(\int\displaylimits_{\{\Phi_2<\infty\}}\abs{\nabla \Phi_2 f}^2
%%		\,d(x,v)\right)^
%%		{\frac{\alpha}{2}}\gamma\left(\int\displaylimits_{\{\Phi_2<\infty\}}
%%		\abs{f}^2\,d(x,v)\right)^{\frac{2-\alpha}{2}}\notag\\
%		&\leq \frac{\alpha}{2\gamma^{\frac{2}{\alpha}}}\int\displaylimits_{\{\Phi_2<\infty\}}
%		\abs{\nabla \Phi_2f}^2\,d(x,v)
%		+\frac{(2-\alpha)\gamma^{\frac{2}{2-\alpha}}}{2}\int\displaylimits_{\{\Phi_2<\infty\}}
%		\abs{f}^2\,d(x,v).
%		\end{align}
%		Since $\gamma$ was choosen arbitrary we can set $\gamma=(2\alpha K)^{\frac{\alpha}{2}}$ 
%		and obtain 
		\begin{equation}\label{alpha}
		(\abs{\nabla_v \Phi_2}^\alpha f,f)\leq
		\frac{1}{4K}\norm{\abs{\nabla_v \Phi_2}f}^2
		+\frac{(2-\alpha)(2\alpha K)^{\frac{\alpha}{2-\alpha}}}{2}\norm{f}^2.
		\end{equation}
		Consequently, for $f\in D_0$ the inequality $(\ref{nablaPhi2})$ becomes 
%		\begin{align}
%		\norm{\abs{\nabla_v\Phi_2} f}^2 &\leq
%		4(-L^\prime_0f,f) + 2K\left(1+\frac{(2-\alpha)(2\alpha K)^{\frac{\alpha}{2-\alpha}}}{2}\right)
%		\norm{f}^2 +\frac{1}{2}\norm{\abs{\nabla_v\Phi_2}f}^2\notag
%		\end{align}
		\begin{align}
		\norm{\abs{\nabla_v\Phi_2} f}^2 \leq
		8(-L^\prime_0f,f) + C\norm{f}^2,\label{Opbou} 
		\end{align}
		with $C=4K(1+\frac{(2-\alpha)(2\alpha K)^{\frac{\alpha}{2-\alpha}}}{2})$. Since 
		(\ref{Opbou}) holds
		we conclude that $\abs{\nabla_v\Phi_2}P_k$ is $L_k$ bounded with $L_k$-bound zero 
		and so is $B_0^k$ for each $k\in\mathbb N$.
		Now we are able to apply \cite[Lemma 3]{CG08}
		implying essential m-dissipativity of 
		\begin{equation}
		\left(L^\prime,D_0\right):=\left(L_0^\prime+B_0,D_0\right)=
		\left(\Delta_v-\frac{\abs{\nabla_v\Phi_2}^2}{4}
		+\frac{\Delta_v\Phi_2}{2}+i\nabla_v\Phi_2\cdot x,D_0\right). 
		\end{equation}
		Note: The estimates (\ref{nablaPhi2}),(\ref{nablaPhizweizwei}),(\ref{alpha}) and (\ref{Opbou}) also hold for $ f $ in the bigger space $ L^2(\mathbb{R}^d)\otimes U_{\Phi_2} C_c^{\infty}(\{\Phi_2<\infty\}) $.
		
%		Next we use the Fourier transform in $x$ direction as a unitary map to transform 
%		the added operator $i\nabla_v\Phi_2\cdot x$ into $\nabla_v\Phi_2\cdot \nabla_x$. Therefore 
		The set $D_1=C^\infty_c(\Rd{d})\otimes U_{\Phi_2} C_c^{\infty}(\{\Phi_2<\infty\})$ (analogue definition as for $ D_0 $)
% But $ D_1 $ 
		forms a core for the closure of $(L^\prime,D_0)$, hence, $(L^\prime,D_1)$ is essentially
		m-dissipative, too.	The extension of $(L^\prime,D_1)$ to $D_2=\mathcal{S}(\mathbb{R}^{d})\otimes
		U_{\Phi_2} C_c^{\infty}(\{\Phi_2<\infty\})$ is still dissipative, hence the closure of $(L^\prime,D_2)$ is a dissipative extension 
		of the closure of $(L^\prime,D_1)$ and therefore their closures coincide by \cite[Chapter 1, Remark 3.8]{G85}, i.e.,
		\begin{equation}\label{SC}
		\left(L^\prime,\mathcal{S}(\mathbb{R}^{d})\otimes
		U_{\Phi_2} C_c^{\infty}(\{\Phi_2<\infty\})\right)
		\text{is essentially m-dissipative.} 
		\end{equation}
		Denote by $\mathcal{F}$ the Fourier transform on $L^2(\Rd{d})$.
		Recall the well-knonwn property of $\mathcal{F}:$
%		 (see e.g. \cite[section 9.3]{Do10}) :
%		\begin{enumerate}
%			\item[($\mathcal F_1$)] $\mathcal F$ is an isometric isomorphism on $L^2(\Rd{d})$,
%			\item[($\mathcal F_2$)] $\mathcal F(\mathcal S(\Rd{d}))=\mathcal S(\Rd{d})$,
%			\item[($\mathcal F_3$)] For 
%			$f\in\mathcal{S}(\mathbb{R}^{d})$ and a multi-index	$\alpha\in\mathbb{N}_0^d$ it holds
			\begin{equation}\label{Fourier}
			\mathcal{F}^{-1}(x^s f)=(-i)^{\abs{s}}\partial^s(\mathcal{F}^{-1}f),\text{ for }f\in\mathcal{S}(\mathbb{R}^{d})\text{ and }s\in\mathbb{N}_0^d.
			\end{equation}
%		\end{enumerate}
		Let $f=f_1\otimes f_2\in D_2$. Define $\mathcal{F}_xf:=\mathcal{F}f_1\otimes f_2$
		and extend $\mathcal{F}_x$ linearly to $D_2$ and afterwards 
%		$\mathcal{F}_x$ is well-defined because 
%		$\mathcal{F}$ is linear and isometric.
%		In particular $\mathcal{F}_x$ is an isometry, too.
%		Since $D_2$ is dense in $L^2(\{\Phi_2<\infty\})$ and is left invariant under $ \mathcal F_x $
%		(see remark on page $\pageref{dense}$)   
		to a unitary transformation on $L^2\left(\{\Phi_2<\infty\}\right)$ (similarly as one does for $ \mathcal F $) which we also denote by $\mathcal{F}_x$. 
		$\mathcal F_x$ leaves the set $D_2$ invariant, because $\mathcal{S}(\Rd{d})$  is invariant under $\mathcal F$.	Using the identity $(\ref{Fourier})$ one obtains
%		 easily checks that $(L^\prime,D_2)$ transforms
%		under $\mathcal{F}_x$ into the essentially m-dissipative operator $( L,D_2)$ given by
		\begin{equation}
		 \tilde Lf=\mathcal F_x^{-1}L^\prime\mathcal F_x f=\left(\Delta_v+\frac{\Delta_v\Phi_2}{2}-\frac{\abs{\nabla_v\Phi_2}^2}{4}
		+\nabla_v\Phi_2\cdot\nabla_x\right)f,\quad f\in D_2. 
		\end{equation}		
		We perturb $ \tilde L $ with the antisymmetric operator $(B_1,D_2)$ given by
		$ B_1f=\sum_{i=1}^d \partial_{x_i}\Phi_1\partial_{v_i}f$, $ f\in D_2$. Since $ \Phi_1 $ is Lipschitz continuous $(B_1,D_2)$ is well-defined.
%		, i.e., $B_1f\in L^2(\{\Phi_2<\infty\})$ for $f\in D_2$.
%		Furthermore it holds $\abs{\nabla_x\Phi_1}\leq C_{\Phi_1}$
%		(see e.g. \cite[Satz 5.23]{Do10}). 
%		Since $\Phi_2(x,v)=\Phi_2(v)$ we have $\nabla_x\Phi_2=0$. 
%		Applying integration by parts for Sobolev functions we see that 
%		$\nabla_v\Phi_2\cdot\nabla_x$ is antisymmetric.
%		It follows that $(\nabla_v\Phi_2\cdot\nabla_xf,f)$ is purely
%		imaginary for $f\in D_2$.
		As in the derivation of (\ref{Opbou}) we obtain finite constants $ C_1$ and $C_2 $ such that
%		 As above let $K$ and $\alpha$ be the constants from assumption
%		$(\Phi_25)$ as well as $f\in D_2$. With the estimate (\ref{+Anti}) in mind it holds
		\begin{align}
		\norm{B_1f}=\norm{\nabla_x\Phi_1\cdot\nabla_v f}^2 &\leq C^2_{\Phi_1}
		\sum_{i=1}^d(\partial_{v_i} f,\partial_{v_i} f)=C^2_{\Phi_1}(-\Delta_v f,f)\notag\\
%		&\leq C^2_{\Phi_1}\left((-\Delta_v+\frac{\abs{\nabla_v\Phi_2}^2}{4}-
%		\frac{\Delta_v\Phi_2}{2})f,f\right)+\frac{C^2_{\Phi_1}}{2}(\Delta_v\Phi_2 f,f)\notag\\
%		&\leq C^2_{\Phi_1}\left|\underbrace{\left((-\Delta_v+\frac{\abs{\nabla_v\Phi_2}^2}{4}-
%			\frac{\Delta_v\Phi_2}{2})f,f\right)}_{\in \mathbb R}
%		-\underbrace{\left(\nabla_v\Phi_2\cdot\nabla_xf,f\right)}_{\in i\mathbb R}\right|\notag\\
%		&\hspace{4.5mm}+\frac{C^2_{\Phi_1}}{2}(\Delta_v\Phi_2 f,f)\notag\\
%%		&=\leq C^2_{\Phi_1}\abs{(- L f,f)}+
%		\frac{C^2_{\Phi_1}}{2}(\Delta_v\Phi_2 f,f)\notag\\
		&\leq C_1(-L_0^\prime f,f)+C_2\norm{f}^2.\label{Lastesti}
		\end{align}
%		The term $(\abs{\nabla_v\Phi_2}^\alpha f,f)$ was already estimated above in (\ref{alpha})
%		but only for $f\in D_0$. 
%		Anyway, the arguments used work for $f\in D_2$, too. We obtain 
%		\begin{align}
%		\norm{\nabla_x\Phi_1\cdot\nabla_v f}^2 &\leq C^2_{\Phi_1}\abs{(- L f,f)}+
%		\frac{C^2_{\Phi_1}K}{2}\left(1+\frac{(2-\alpha)(2\alpha K)^{\frac{2}{2-\alpha}}}{2}\right)
%		\norm{f}^2\notag\\
%		&+\frac{1}{4K}\norm{\abs{\nabla \Phi_2}f}^2.
%		\end{align}
%		Let us make a fruitful observation. 
%		concerning the inequality (\ref{Opbou}).
%		The same arguments we used for establishing (\ref{Opbou}) work for
%		$f\in D_2$ instead of $D_1$, too.
		 Since $(L_0^\prime,D_2)$ is symmetric it holds that 
		$(L_0^\prime f,f)\in \mathbb R$, for $ f\in D_2 $. Let $A$ be an arbitrary antisymmetric linear operator on $D_2$. In particular, for $f\in D_2$ it holds that $(Af,f)\in i\mathbb R$. Hence one obtains
		\begin{equation}\label{+Anti}
		(-L^\prime_0f,f)\leq\abs{\underbrace{(-L^\prime_0f,f)}_{\in\mathbb R}+\underbrace{(Af,f)}_{\in i\mathbb R}}.
		\end{equation}	
		Applying the inequality $(\ref{+Anti})$ for the choice $A=-\nabla_v\Phi_2\cdot\nabla_x$ to (\ref{Lastesti}) one concludes
%		\begin{equation}
%		\norm{\abs{\nabla \Phi_2}f}^2\leq 8\abs{(- Lf,f)}+C\norm{f}^2. 
%		\end{equation}
%		We finally obtain
		\begin{equation}
		\norm{\nabla_x\Phi_1\cdot\nabla_v f}^2 
		\leq C_1\abs{(- \tilde L f,f)}+C_2\norm{f}^2.
		\end{equation}
		By \cite[Chapter 3.1, Lemma 3.9]{D8} we deduce that 
		\[L=\tilde L-\nabla_x\Phi_1\cdot\nabla_v=\Delta_v-\frac{\abs{\nabla_v\Phi_2}^2}{4}+
		\frac{\Delta_v\Phi_2}{2}+\nabla_v\Phi_2\cdot\nabla_x-\nabla_x\Phi_1\cdot\nabla_v\]
		defined on $D_2$ is essentially m-dissipative on $L^2(\{\Phi_2<\infty\})$.

		We apply $(\ref{+Anti})$ with 	$ A=-\nabla_v\Phi_2\cdot\nabla_x+\nabla_x\Phi_1\cdot\nabla_v$ to extend (\ref{Opbou}) for $ L $ instead of $ L_0^\prime $, i.e., 
		% that 
		% $(\ref{Opbou})$
		% also holds for $L$ instead of $ L_0^\prime$, i.e.
		\begin{equation}\label{Opbou3}
		\norm{\abs{\nabla\Phi_2}f}^2\leq r\abs{(Lf,f)}+M\norm{f}^2,\quad f\in D_2,
		\end{equation}
		for finite constants $ r,M $.		
		We restrict $ L $ to $ D_1 $ and observe that essential m-dissipativity is preserved, since $ C^\infty_c(\Rd{d}) $ is dense in $ \mathcal S(\Rd{d}) $ (w.r.t. the Schwartz space topology on $ \mathcal S(\Rd{d}) $).		
		Now we transform via the adjoint of unitary map from (\ref{Uni}), i.e., 
		\begin{equation}\label{unitrafo}
		U^*:L^2(\{\Phi_2<\infty\})\longrightarrow L^2(\Rd{2d},\mu_{\mathbf\Phi}),
		f\mapsto e^{\frac{\Phi_1+\Phi_2}{2}}\tilde f,
		\end{equation}
		where $ \tilde f=1_{\{\Phi_2<\infty\}}f$. For $ f=f_1\otimes f_2\in D_1 $ one has $ U^* f= e^{\frac{\Phi_1}{2}}f_1\otimes  e^{\frac{\Phi_2}{2}}f_2 $.
		Denote by $ U^*_{\Phi_1} $ the unitary map 
		$
		U^*_{\Phi_1}:L^2(\Rd{d})\longrightarrow L^2(\Rd{d},\mu_{\Phi_1}), f\mapsto e^{\frac{\Phi_1}{2}}f.
		$
		 Due to (\ref{L1}), (\ref{L0}), the product rule for Sobolev functions and Proposition \ref{H1infty}, it holds that $ U^* $ transforms $ L $ back into $ L_{\mathbf{\Phi}} $, i.e., we obtain the essentially m-dissipative operator
		\begin{equation}\label{HC}
		(U^*LU,U^*D_1)=\left(L_{\mathbf{\Phi}},
		U^*_{ \Phi_1}C^{\infty}_c(\mathbb{R}^d)\otimes C_c^{\infty}(\{\Phi_2<\infty\})\right). 
		\end{equation}
%		 , i.e, elements vanishing outside some compact set. 
		For $ f\in U^*_{ \Phi_1}C^{\infty}_c(\mathbb{R}^d)\otimes C_c^{\infty}(\{\Phi_2<\infty\}) $ it holds $ Uf\in D_1 $ and hence through (\ref{Opbou3}) we obtain 
		\begin{align}
		\norm{\abs{\nabla_v\Phi_2}f}^2_{\mu_{\mathbf\Phi}}=\norm{\abs{\nabla_v\Phi_2}Uf}^2 
		&\leq r\abs{(- LUf,Uf)}+M\norm{Uf}^2\notag\\
		&=r\abs{(- L_{\mathbf{\Phi}}f,f)_{\mu_{\mathbf\Phi}}}+M\norm{f}_{\mu_{\mathbf\Phi}}^2.\label{Kato2}
		\end{align}
		The lemma of Fatou guarantees that $ (\ref{Kato2}) $ also holds for $ f $ from the closure of $ (\ref{HC}) $.
		To finish the first part we show that $ C^{\infty}_c(\mathbb{R}^d)\otimes C_c^{\infty}(\{\Phi_2<\infty\}) $ is a domain of essential m-dissipativity for $ L_{\mathbf{\Phi}} $. Since $  (L_{\mathbf{\Phi}},C^{\infty}_c(\mathbb{R}^d)\otimes C_c^{\infty}(\{\Phi_2<\infty\})) $ is dissipative by Corollary \ref{invariant} it suffices due to the essential m-dissipativity of $ (\ref{HC}) $ and \cite[Chapter 1, Remark 3.8]{G85} to show that the closure of $  (L_{\mathbf{\Phi}},C^{\infty}_c(\mathbb{R}^d)\otimes C_c^{\infty}(\{\Phi_2<\infty\})) $ is an extension of $ (\ref{HC}) $.
		To this end let $ f=f^1\otimes f^2\in 	U^*_{ \Phi_1}C^{\infty}_c(\mathbb{R}^d)\otimes C_c^{\infty}(\{\Phi_2<\infty\}) $. Observe that $ U^*_{ \Phi_1}C^{\infty}_c(\mathbb{R}^d)
		 $ is by Proposition \ref{H1infty} a subset of $ H^{1,2}(\Rd{d})$. Choose a sequence $ (f_n^1)_{n\in\mathbb N}$ from $ C^\infty_c(\Rd{d})$ such that $ {f_n^1\longrightarrow f^1} $
		in $H^{1,2}(\Rd{d})$ and $\text{supp}(f^1_n)\subseteq K$, $K\subseteq \Rd{d}$ compact and independent of $ n $ which is possible since $ f^1 $ is already compactly supported. For $f_n:=f_n^1\otimes f^2$, $ n\in\mathbb N $, it holds by construction and the fact that the density $ e^{-\Phi_1-\Phi_2} $ of $ \mu_{\mathbf\Phi} $ is locally bounded that $f_n\longrightarrow f $, $ L_{\mathbf{\Phi}}f_n\longrightarrow L_{\mathbf{\Phi}}f $ and $ \abs{\nabla_v\Phi_2}f_n\longrightarrow \abs{\nabla_v\Phi_2}f $ in $ \mathcal H_{\mathbf{\Phi}} $ as $ n\to\infty $. This shows that $C^\infty_c(\Rd{d})\otimes C_c^{\infty}(\{\Phi_2<\infty\})$ is a core for the closure of (\ref{HC}). 

	\noindent\textbf{2nd part:}\\
		Let $\Phi_1$ be locally Lipschitz continuous. Dissipativity is due to Corollary \ref{invariant}. 
%		that $\left(L_{\mathbf{\Phi}},C_c^{\infty}(\{\Phi_2<\infty\})\right)$ is dissipative.
		To prove m-dissipativity we show that  $(1-L_{\mathbf{\Phi}})C_c^{\infty}(\{\Phi_2<\infty\})$ is dense.
		Since $C_c^{\infty}(\{\Phi_2<\infty\})$ is dense it suffices to approximate
		$0\neq g\in C_c^{\infty}(\{\Phi_2<\infty\})$. Let  
		$f \in C_c^\infty(\{ \Phi_2 <\infty \})$ be arbitrary and $\epsilon>0$. By the compactness of
		the support of $g$ we can choose cut off functions $\chi,\nu\in
		C_c^{\infty}(\mathbb{R}^d)$ %\{\Phi_2<\infty\}\subseteq.
		such that the functions defined by $\chi(x,v)=\chi(x)$, $\nu(x,v)=\nu(x)$ fulfil the properties $0\leq\chi\leq\nu\leq 1$, $\chi\equiv1$ on $\text{supp}(g)$, $ \nu\equiv1$ on $\text{supp}(\chi)$.
		It holds that $L_{\mathbf{\Phi}}(\chi f)=
		\chi L_{\mathbf{\Phi}}f+f\nabla_v \Phi_2\cdot\nabla_x\chi$ since $\nabla_v\chi=0$.
		By the choice of $\nu$ and $\chi$ we obtain
%		it holds
%		\begin{equation}\label{tildenu}
%		\nu(x)\Phi_1(x)=\Phi_1(x),\quad \forall x\in \text{supp}(\chi)
%		\end{equation}
%		and consequently
%		\begin{equation}\label{nu}
%		(\nu\Phi_1)(x,v)=\Phi_1(x,v),\quad \forall(x,v)\in\text{supp}(\chi).  
%		\end{equation}
%		Furthermore it holds that $L_{\mathbf{\Phi}}(\chi f)=
%		\chi L_{\mathbf{\Phi}}f+f\nabla_v \Phi_2\cdot\nabla_x\chi$ since $\nabla_v\chi=0$. 
		\begin{equation}\label{Kato3}
		\norm{(1-L_{\mathbf{\Phi}})(\chi f)-g}_{\mu_{\mathbf\Phi}}
%		&\leq
%		\norm{\chi(1-L_{\mathbf{\Phi}})f-\chi g}_{\mu_{\mathbf\Phi}}
%		+\norm{f\nabla_v\Phi_2\cdot\nabla_x\chi}_{\mu_{\mathbf\Phi}}\notag \\ 
%		(\ref{nu})\longrightarrow&= \norm{\chi\left((1-L_{\nu\Phi_1+\Phi_2})f- g\right)}_{\mu_{\nu\Phi_1+\Phi_2}}
%		+\norm{f\nabla_v\Phi_2\cdot\nabla_x\chi}_{\mu_{\nu\Phi_1+\Phi_2}}\notag\\ 
%		&\leq \norm{(1-L_{\nu\Phi_1+\Phi_2})f-g}_{\mu_{\nu\Phi_1+\Phi_2}}\notag\\ 
%		&\hspace{4.5mm}+\sum_{i=1}^d\norm{\partial_i\chi}_\infty\norm{f\partial_{v_i}\Phi_2}
%		_{\mu_{\nu\Phi_1+\Phi_2}}\notag\\ 
		\leq \norm{(1-L_{(\nu\Phi_1,\Phi_2)})f-g}_{\mu_{(\nu\Phi_1+\Phi_2)}}
		+\norm{f\abs{\nabla_v\Phi_2}}
		_{\mu_{(\nu\Phi_1+\Phi_2)}}\sum_{i=1}^d\norm{\partial_i\chi}_\infty
		\end{equation}
		Since $\nu\Phi_1$ is globally Lipschitz continuous
%		. We make again use of the
%		characterization of $H^{1,\infty}$ as the the set of equivalence classes having Lipschitz continuous representitives (see Remark \ref{densedefined}(i)).
%		Since $\Phi_1$ is locally Lipschitz continuous and therefore an element of 
%		$H^{1,\infty}_{loc}(\Rd{d})$ 
%		the product $\nu\Phi_1$ is in $H^{1,\infty}(\Rd{d})$.
%		Hence there is a Lipschitz continuous function which coincides with 
%		$\nu\Phi_1$ a.e., but $\nu\Phi_1$ is continuous and 
%		so we can conclude that $\nu\Phi_1$ is indeed Lipschitz continuous.
		we can use the first part and therein the inequality $(\ref{Kato2})$ to estimate the last term in $(\ref{Kato3})$ by
		\begin{align}
		\norm{f\abs{\nabla\Phi_2}}_{\mu_{(\nu\Phi_1,\Phi_2)}}
%		&\leq\left(A\norm{-L_{\nu\Phi_1+\Phi_2}f}
%		_{\mu_{\nu\Phi_1+\Phi_2}}
%		+B\norm{f}_{\mu_{\nu\Phi_1+\Phi_2}}\right)\notag\\
%		&\leq\left(A\norm{(1-L_{\nu\Phi_1+\Phi_2})f}
%		_{\mu_{\nu\Phi_1+\Phi_2}}+(B+A)\norm{f}_{\mu_{\nu\Phi_1+\Phi_2}}\right)\notag\\
		&\leq C\left(\norm{(1-L_{(\nu\Phi_1,\Phi_2)})f}_{\mu_{(\nu\Phi_1+\Phi_2)}}+
		\norm{f}_{\mu_{\nu\Phi_1+\Phi_2}}\right) \tag{2.29a}\label{Final}
		\end{align}
		for some positive, finite constant $C$. Since $L_{\nu\Phi_1,\Phi_2}$ is dissipative it holds
%		Recall the decomposition of $L_{\nu\Phi_1+\Phi_2}$ from Lemma \ref{SymAntisy} into
%		$L_{\nu\Phi_1+\Phi_2}=S+A$, with negative definite $S$ and antisymmetric $A$.
		\begin{equation} 
		(f,f)_{\mu_{(\nu\Phi_1+\Phi_2)}}\leq\abs{((1-L_{(\nu\Phi_1,\Phi_2)})f,f)_{\mu_{(\nu\Phi_1+\Phi_2)}}}\leq\norm{f}_{\mu_{(\nu\Phi_1+\Phi_2)}}\norm{(1-L_{(\nu\Phi_1,\Phi_2)})f}_{\mu_{(\nu\Phi_1+\Phi_2)}}\tag{2.29b}.
		\end{equation}
		\stepcounter{equation}
		Now,  (2.29a)  and  (2.29b) imply 
		\begin{equation}
		\norm{f\abs{\nabla\Phi_2}}_{\mu_{(\nu\Phi_1+\Phi_2)}}\sum_{i=1}^d\norm{\partial_i\chi}_\infty
		\leq  2C
		\left(\norm{(1-L_{(\nu\Phi_1,\Phi_2)})f}_{\mu_{(\nu\Phi_1+\Phi_2)}}\right)\sum_{i=1}^d\norm{\partial_i\chi}_\infty.
		\end{equation}
		The inequality (\ref{Kato3}) becomes 
		\begin{align*}
		\norm{(1-L_{\mathbf{\Phi}})(\chi f)-g}_{\mu_{\mathbf\Phi}}&\leq 
		\norm{(1-L_{(\nu\Phi_1,\Phi_2)})f-g}_{\mu_{(\nu\Phi_1+\Phi_2)}}\\
		&\hspace{4.5mm}+2C\left (\norm{(1-L_{(\nu\Phi_1,\Phi_2)})f}_{\mu_{(\nu\Phi_+\Phi_2)}}\right )
		\sum_{i=1}^d\norm{\partial_i\chi}_\infty\\
		&\leq \norm{(1-L_{(\nu\Phi_1,\Phi_2)})f-g}_{\mu_{(\nu\Phi_1+\Phi_2)}}\\
		&\hspace{4.5mm}+2C(\norm{(1-L_{(\nu\Phi_1,\Phi_2)})f-g}_{\mu_{(\nu\Phi_1+\Phi_2)}}
		+\norm{g}_{\mu_{(\nu\Phi_1+\Phi_2)}})
		\sum_{i=1}^d\norm{\partial_i\chi}_\infty
		\end{align*}
		
		Now we specify our choice of $\chi$. Let $\chi$ be chosen in such a way that
		$\sum_{i=1}^d\norm{\partial_i\chi}_\infty\leq\frac{\epsilon}{8C\norm{g}_{\mu_{\mathbf\Phi}}}$.
		Now $\chi,\nu$ are fixed. By the first part of the proof we know that $L_{\nu\Phi_1+\Phi_2}$ is
		essentially m-dissipative. Therefore
		we can choose an element $f\in C_c^{\infty}(\{\Phi_2<\infty\})$ such that
		$\norm{(1-L_{\nu\Phi_1,\Phi_2})f-g}_{\mu_{\nu\Phi_1+\Phi_2}}<\inf\{\frac{\epsilon}{2}
		,\norm{g}_{\mu_{\nu\Phi_1+\Phi_2}}\}$ and we finally obtain
		\begin{equation*}
		\norm{(1-L_{\mathbf{\Phi}})(\chi f)-g}_{\mu_{\mathbf\Phi}}<\epsilon.
		\end{equation*}
		So far we showed that the closure $ (L_{\mathbf{\Phi}},D(L_{\mathbf{\Phi}})) $ of $(L_{\mathbf{\Phi}},C^\infty_c(\{\Phi_2<\infty\})$
		is the generator of a strongly continuous semigroup of contractions $\left(T_t^{\mathbf{\Phi}}\right)_{t\geq0}$. The Dirichlet property (see \cite[Definition I.4.1]{MaR92} for the definition) of
		$ (L_{\mathbf{\Phi}},D(L_{\mathbf{\Phi}})) $ follows by
		\cite[Lemma 1.9, App. B]{Eb99} and hence by \cite[Proposition I.4.3]{MaR92} the semigroup $\left(T_t^{\mathbf{\Phi}}\right)_{t\geq0}$ is sub-Markovian.
	\end{proof}
\end{T}
\begin{Rem}
From the proof of Theorem \ref{maintheorem} one sees that the condition $ (\Phi_25) $ can also be extended to $ \alpha=2 $ and $ 0\leq K< \frac{1}{2} $.
\end{Rem}
Recalling the decomposition 
from Proposition \ref{SymAntisy} we obtain that for the adjoint $ \left(\hat L_{\mathbf{\Phi}},D(\hat L_{\mathbf{\Phi}})\right) $ of $ \left(L_{\mathbf{\Phi}},D(L_{\mathbf{\Phi}})\right) $ it holds 
\begin{equation}\label{adjdec}
C^\infty_c(\{\Phi_2<\infty\})\subseteq D(\hat L_{\mathbf{\Phi}}),\quad \hat L_{\mathbf{\Phi}}f=Sf-Af,\quad  f\in  C^\infty_c(\{\Phi_2<\infty\}).
\end{equation}
For a symmetric velocity potential $  \Phi_2 $, i.e., $ \Phi_2(v)=\Phi_2(-v),\forall v\in\Rd{d}$,  we can use the velocity reversal as in \cite[p. 153]{C10} , i.e., the unitary transformation on $\mathcal H_{\mathbf{\Phi}} $ given by 
\begin{equation}\label{velorev}
U:\mathcal H_{\mathbf{\Phi}}\longrightarrow\mathcal H_{\mathbf{\Phi}},[f]\mapsto [(x,v)\mapsto f(x,-v)]
\end{equation}
to transform $(L_{\mathbf{\Phi}},C^\infty_c(\{\Phi_2<\infty\}) $ into the operator $ (UL_{\mathbf{\Phi}}U, UC^\infty_c(\{\Phi_2<\infty\})=\left(\hat L_{\mathbf{\Phi}},C^\infty_c(\{\Phi_2<\infty\})\right)$. This implies that the latter is also an essential m-dissipative operator. Hence, the closure of $ \left(\hat L_{\mathbf{\Phi}},C^\infty_c(\{\Phi_2<\infty\})\right) $ coincides with the adjoint of the closure of $(L_{\mathbf{\Phi}},C^\infty_c(\{\Phi_2<\infty\}) $.
% To guarantee that the adjoint  $(\hat L_{\mathbf{\Phi}},D(\hat L_{\mathbf{\Phi}})) $ really coincides with the closure of $\left ( \hat L_{\mathbf{\Phi}}, C^\infty_c(\{\Phi_2<\infty\})\right ) $ one can employ the adjoint semigroup $ \left(\hat T^{\Phi_1,\Phi_2}_t\right)_{t\geq 0}$. Indeed, due to \cite[Lemma 1.1.25 (iii)]{C10} the weak generator of $ \left(\hat T^{\mathbf{\Phi}}_t\right)_{t\geq 0}$ is the adjoint of 
% $(L_{\mathbf{\Phi}},D(L_{\mathbf{\Phi}})) $ and is therefore dissipative (see \cite[proof of Lemma 1.1.23]{C10}), hence, the closure of  $( \hat L_{\mathbf{\Phi}}, C^\infty_c(\Rd{2d}) $ and $(\hat L_{\mathbf{\Phi}},D(\hat L_{\mathbf{\Phi}})) $ are the same.
Therefore, we assume in the following the additional assumption:
\begin{C}\label{C3}
	\mbox{}\vspace{-\topskip}
	\begin{itemize}[align=parleft, labelsep=0.5cm] 				
		\item [$(\Phi_2\thePHIzwei)$]\tabto{1cm}
		$\Phi_2$ is symmetric, i.e., $ \Phi_2(v)=\Phi_2(-v)$, for all $v\in\Rd{d}$. \label{symmetry}
		\stepcounter{PHIzwei}
	\end{itemize}
\end{C}
The next corollary recaps the previous discussion.
\begin{Co}\label{Adjoint}
	Under the assumptions of Theorem \ref{maintheorem} and the additional assumption $ (\Phi_26) $ the formal adjoint $(\hat L_{\mathbf{\Phi}}, C^\infty_c(\{\Phi_2<\infty\})) $ is also an essentially m-dissipative Dirichlet operator. Furthermore, its closure coincides with the adjoint of $(L_{\mathbf{\Phi}},D(L_{\mathbf{\Phi}})) $.
\end{Co}

\subsection{M-Dissipativity for singular $ \Phi_1 $ on $  L^1(\Rd{2d},\mu_{\mathbf\Phi}) $}
In this part we merely assume $ (\Phi_12)-(\Phi_14)^q $, $ q\in[2,\infty] $, for $ \Phi_1 $ and $ (\Phi_21)-(\Phi_26) $ for $ \Phi_2 $. Observe that due to Corollary \ref{invariant} the operator $ \left(L_{\mathbf{\Phi}},C^\infty_c(\{\Phi_1,\Phi_2 <\infty\})\right) $ is closable on $ L^1(\Rd{2d},\mu_{ \mathbf\Phi}) $ and its closure $ \left(L_{\mathbf{\Phi},1},D(L_{\mathbf{\Phi},1})\right) $ is dissipative.
% Observe that the upper index $ (1) $ only indicates that we consider this operator on the space $ L^1(\Rd{2d},\mu_{ \mathbf\Phi}) $. Later we use the same notation to describe the case when we consider this operator on the spaces $ L^p(\Rd{2d},\mu_{ \mathbf\Phi})  $ for $ p\geq1 $, which is different to the notation used in (\ref{Lbg}) where the upper index is without brackets.  
%is to approximate the potential $ \Phi_1 $ by locally Lipschitz continuous functions $ \Phi_{1,n} $ and use the results from the previous part. To apply the technique from  we have to provide some prerequisites.
%
The next proposition is taken from \cite[Lemma 3.7]{CG10}. We only state the parts which are necessary for our needs.
\begin{Prop}\label{extensionweak}
	The set $ C^\infty_c(\{\Phi_2 <\infty\}) $ is contained in $ D(L_{\mathbf{\Phi},1}) $ and for $ f\in C^\infty_c(\{\Phi_2 <\infty\}) $ it holds $ L_{\mathbf{\Phi},1}f=L_{\mathbf{\Phi}}f $.
\end{Prop}

%\begin{Rem}
%The proof of Proposition \ref{extensionweak} is the only time we need the assumption $ (\Phi_13) $.
%\end{Rem}

\begin{Co}\label{prepmdis}
	$ \left(L_{\mathbf{\Phi}},C^\infty_c(\{\Phi_1,\Phi_2 <\infty\})\right) $ is essentially m-dissipative on $ L^1(\Rd{2d},\mu_{ \mathbf\Phi}) $ iff its extension $ \left(L_{\mathbf{\Phi}},C^\infty_c(\{\Phi_2 <\infty\})\right) $ is.
\end{Co}
\noindent The next lemma provides a sequence of smooth potentials $ (\Phi_{1,n})_{n\in\mathbb N} $ approximating $ \Phi_1 $
in a suitable sense. See \cite[Lemma 3.10]{CG10} for the proof.
\begin{Le}\label{Phin}
	Let $ \Phi=\Phi_1 $ fulfill $ (\Phi_12) $, $ (\Phi_13) $, $ (\Phi_14)^q $. Then there exist smooth $ \Phi_{n}=\Phi_{1,n} $ such that $ \Phi_{n}\leq \Phi $ and $ \nabla\Phi_{n} \overset{n\to\infty}{\longrightarrow}\nabla\Phi$ in $ L^q_{loc}(\Rd{d},\mu_{\Phi}) $. Furthermore, the family $ (\Phi_n)_{n\in\mathbb N} $ is uniformly bounded from below.
\end{Le}

\noindent In the following we assume additionally on $ \Phi_2 $:
\begin{C}
	\mbox{}\vspace{-\topskip}
	\begin{itemize}[align=parleft, labelsep=0.5cm]		
		\item[$(\Phi_2\thePHIzwei)$]\tabto{1cm}
		\stepcounter{PHIzwei}
		$\mu_{{\Phi}_2}$ is a finite measure, i.e., 
		$\mu_{{\Phi}_2}(\Rd{d})=
			\int\displaylimits_{\Rd{d}}e^{-\Phi_2}\,dv<\infty.  
			$
		\item[$(\Phi_2\thePHIzwei)$]\tabto{1cm}
			The measurable function $\abs{\nabla\Phi_2}$ is square integrable w.r.t.
			$\mu_{{\Phi}_2}$, i.e.,\tabto{1cm}
			{$
				\int_{\Rd{d}}\abs{\nabla\Phi_2}^2\,d\mu_{{\Phi}_2}=
				\int_{\Rd{d}}\abs{\nabla\Phi_2}^2e^{-\Phi_2}\,dv<\infty.
				$}
			\stepcounter{PHIzwei} 
	\end{itemize}
\end{C}

\begin{T}\label{mdissiL1}
	Assume $ (\Phi_12)-(\Phi_14)^q $ and $ (\Phi_21)-(\Phi_28) $. Additionally one of the following assumptions are assumed.
	\begin{enumerate}
		\item $ \mu_{ \mathbf\Phi} $ is a finite measure.
		\item $ (\Phi_14)^q $ holds for $ q>d $.
	\end{enumerate}
%	Then the operator $\left(L_{\mathbf{\Phi}},C_c^{\infty}(\{\Phi_1,\Phi_2<\infty\})\right)$ is essentially m-dissipative on $L^1(\Rd{2d},\mu_{ \mathbf\Phi}) $. The strongly continuous contraction semigroup $ \left(T^{\mathbf\Phi}_t\right)_{t\geq 0}$ generated by the closure
%		$\left(L_{\mathbf{\Phi},1},D(L_{\mathbf{\Phi},1})\right)$
%		is sub-Markovian. 
	Then the operator $\left(L_{\mathbf{\Phi},1},D(L_{\mathbf{\Phi},1})\right)$ generates a strongly continuous contraction semigroup $ \left(T^{\mathbf\Phi}_{t,1}\right)_{t\geq 0}$ on $L^1(\Rd{2d},\mu_{ \mathbf\Phi}) $. Furthermore, this semigroup is sub-Markovian. 
		\begin{proof} 
		Together with Theorem \ref{maintheorem}, Corollary \ref{prepmdis} and Lemma \ref{Phin} we provided all prerequisites to apply the proof of \cite[Theorem 3.11]{CG10}. The sub-Markovian property of $ \left(T^{\mathbf\Phi}_{t,1}\right)_{t\geq 0}$ holds due to \cite[Appendix B, Lemma 1.9]{Eb99}.
		\end{proof}
\end{T}
%The next lemma is taken from \cite{CG10}[Lemma 3.14]
%\begin{Le}\label{Approximation}
%Let $ \Phi_{1,n} $, $ n\in\mathbb N $, be defined as in Lemma \ref{Phin} and denote by $ \mathbf{\Phi}_n=(\Phi_{1,n},\Phi_2) $. For $ f\in L^1(\Rd{2d},d(x,v)) $ and $ t\geq0 $ it holds
%\begin{equation}\label{semigroupap}
%\norm{T^{\mathbf{\Phi}}_tf-T^{\mathbf{\Phi}_n}_tf}_{L^1(\Rd{2d},\mu_{ \mathbf\Phi})}\overset{n\to\infty}{\longrightarrow}0.
%\end{equation}
%\end{Le}

\noindent Observe that the velocity reversal $ U $ from  (\ref{velorev}) is also a bijective isometry on the space $L^1(\Rd{2d},\mu_{ \mathbf\Phi}) $. Hence, the closure of the formal adjoint $\left(\hat L_{\mathbf{\Phi}},C_c^{\infty}(\{\Phi_1,\Phi_2<\infty\})\right)$ in $L^1(\Rd{2d},\mu_{ \mathbf\Phi}) $ is the generator of a sub-Markovian stongly continuous contraction semigroup $ \asg{\mathbf{\Phi}}{,1} $ on $L^1(\Rd{2d},\mu_{ \mathbf\Phi}) $. The two semigroups $ \sg{\mathbf{\Phi}}{,1} $ and $ \asg{\mathbf{\Phi}}{,1} $ give rise to contraction semigroups $ \sg{\mathbf{\Phi}}{,p} $ and $ \asg{\mathbf{\Phi}}{,p} $ on $ L^p(\Rd{2d},\mu_{\mathbf\Phi}) $ for every $ p\in[1,\infty] $ which are also strongly continuous for $ p\in[1,\infty) $. These semigroups coincide with $ \sg{\mathbf{\Phi}}{,1} $ and $ \asg{\mathbf{\Phi}}{,1} $ on $L^1(\Rd{2d},\mu_{\mathbf\Phi})\cap L^\infty(\Rd{2d},\mu_{\mathbf\Phi}) $, respectively (see \cite[Lemma 1.3.11]{C10} for details).
%\pagebreak
\begin{Le}\label{Lpgen} Let the assumptions of Theorem \ref{mdissiL1} hold true. Furthermore, let $ p\in[1,\infty) $.
\begin{enumerate}
	\item[(i)]
	The generator $\left(L_{\mathbf{\Phi},p},D(L_{\mathbf{\Phi},p})\right)$ of $  \sg{\mathbf{\Phi}}{,p} $ is given by the closure of $ (L_{\mathbf{\Phi},1},D(L_{\mathbf{\Phi}})_p) $ in $ L^p(\Rd{2d},\mu_{ \mathbf\Phi}) $, where $ D(L_{\mathbf{\Phi}})_p=\left\{f\in D(L_{\mathbf{\Phi},1})\mid f,L_{\mathbf{\Phi},1}f\in L^p(\Rd{2d},\mu_{ \mathbf\Phi})\right\} $. In particular, for $ f\in D(L_{\mathbf{\Phi}})_p$ it holds $ L_{\mathbf{\Phi},p}f=L_{\mathbf{\Phi},1}f $.
	\item[(ii)] The contraction semigroups $ \sg{\mathbf{\Phi}}{,\frac{p}{p-1}} $ and $ \asg{\mathbf{\Phi}}{,\frac{p}{p-1}} $ are the adjoints of $ \asg{\mathbf{\Phi}}{,p} $ and $ \sg{\mathbf{\Phi}}{,p} $, respectively.
	\item[(iii)] The semigroup $ \sg{\mathbf{\Phi}}{,\infty} $ is conservative and $ \mu_{ \mathbf\Phi} $ is invariant for $ \sg{\mathbf{\Phi}}{,1} $, i.e.,
	$ T^{ \mathbf\Phi}_{t,\infty}1=1 $ for all $ t\geq0 $ and $ \int_{\Rd{2d}}T_{t,1}^{\mathbf{\Phi}}f\,d\mu_\mathbf{\Phi}=\int_{\Rd{2d}}f\,d\mu_{\mathbf{\Phi}},\forall f\in L^1(\Rd{2d},\mu_\mathbf{\Phi}),t\geq 0$. The same statements also hold for $ \asg{\mathbf{\Phi}}{,\infty} $ and $ \asg{\mathbf{\Phi}}{,1} $, respectively.
%	
%	\int_{E}S_tf\,d\mu=\int_{E}f\,d\mu,\forall f\in L^1(E,\mu),t\geq 0\Longleftrightarrow			
%	\int_EAf\,d\mu=0, \forall f\in\mathcal C.
\begin{proof}
	For part (i) see \cite[Lemma 1.3.11]{C10}, (ii) works analog as in \cite[Lemma 3.16]{CG10}. We prove part (iii): The invariance of $ \mu_{ \mathbf\Phi} $ for $ \sg{\mathbf \Phi}{,1} $ holds by Corollary \ref{invariant}, i.e., $ \int_{\Rd{2d}}L_{\mathbf{\Phi},1}f\,d\mu_\mathbf{\Phi}=0$, for all $ f\in D(L_{\mathbf{\Phi},1}) $. The same argument proves invariance of $ \mu_{ \mathbf\Phi} $ for $ \asg{\mathbf{\Phi}}{,1} $. The conservativeness follows by (ii) and the invariance of $ \mu_{ \mathbf\Phi} $ for $ \asg{\mathbf{\Phi}}{,1} $ and $  \sg{\mathbf{\Phi}}{,1}  $. 
\end{proof}
\end{enumerate}
\end{Le}

% !TeX spellcheck = en_US
\section[Existence of Martingale solutions for $\left(L_{\mathbf{\Phi},2},D(L_{\mathbf{\Phi},2})\right)$]{Existence of Martingale solutions for $\left(L_{\mathbf{\Phi},2},D(L_{\mathbf{\Phi},2})\right)$}\label{ExPr}
In this section we use the results of \cite[Section 3.4]{CG10} to state the existence martingale solutions for
operator $\left(L_{\mathbf{\Phi},2},D(L_{\mathbf{\Phi},2})\right)$, see Theorem \ref{MP} for the precise statement. The core is the result \cite[Theorem 1.1]{BBR06} which provides a $ \mu_\mathbf{\Phi}- $standard right process which is associated in the resolvent sense with $\left(L_{\mathbf{\Phi},1},D(L_{\mathbf{\Phi},1})\right)$, see also the last mentioned reference for the definition of a $ \mu_\mathbf{\Phi}- $standard right process. Theorem \ref{MP} isn't stated in its full generality as in \cite[Theorem 3.1.(iii)]{CG10}. We restrict ourselves to the cases necessary for the applications in mind from section \ref{WCoL}. The proof is completely analog to the one in \cite{CG10} and is therefore omitted.

\noindent Throughout this paper the spaces of continuous functions $ C\left([0,T],E\right) $, $ C\left([0,\infty),E\right) $, where $ (E,m) $ is a metric space and $ T\in\mathbb N $, are always equipped with the topologies of uniform convergence on compact sets and the respective Borel $ \sigma- $algebras. 

%See \cite[p.642]{CG10} for a discussion on how the assumptions of \cite[Theorem 1.1]{BBR06} can be checked.
%The following theorem can be proven completely analogously as in \cite[Section 3.4]{CG10}. Throughout this paper the spaces of continuous functions $ C\left([0,T],E\right) $, $ C\left([0,\infty),E\right) $, where $ (E,m) $ is a metric space and $ T\in\mathbb N $, are always equipped with the metric
%\begin{align}
%d_T(f,g)&=\sup_{t\in[0,T]}m(f(t),g(t)),\quad  &f,g\in C\left([0,T],E\right),\label{metricinfty}\\  
%d(f,g)&=\sum_{T=1}^\infty2^{-T}\frac{d_T(f_{|[0,T]},g_{|[0,T]})}{1+d_T(f_{|[0,T]},g_{|[0,T]})} ,  &f,g\in C\left([0,\infty),E\right)\label{metricfin}. 
%\end{align}
%Furthermore, they always carry the respective Borel $ \sigma- $algebras.
\begin{T}\label{MP}
	Assume $ (\Phi_12)-(\Phi_14)^2,(\Phi_15),(\Phi_16) $, $ (\Phi_21)-(\Phi_28) $.
	Let $ 0\leq h\in L^1(\Rd{2d},\mu_{ \mathbf\Phi})\cap L^2(\Rd{2d},\mu_{ \mathbf\Phi}) $ be a probability density w.r.t. $ \mu_{ \mathbf\Phi} $. Denote by $ \langle\cdot,\cdot \rangle_{\mu_{\mathbf{\Phi}}} $ the dual pairing between $ L^1(\Rd{2d},\mu_{ \mathbf\Phi})$ and $ L^\infty(\Rd{2d},\mu_{ \mathbf\Phi})$.
	There exists a probability law $\mathbb P_{h\mu_{\mathbf{\Phi}}} $ with initial distribution $ h\mu_{\mathbf{\Phi}} $
	on $C([0,\infty),\{\Phi_1,\Phi_2<\infty\}) $ which is associated with the semigroup $ \left(T^{\mathbf{\Phi}}_{t,1}\right)_{t\geq 0}$, i.e., for all
	$ f_1,...,f_k\in L^\infty(\Rd{2d},\mu_{ \mathbf\Phi}) $ and $ 0\leq t_1<...< t_k $, $ k\in\mathbb N $, it holds
	\begin{equation}\label{associated}
	\mathbb{E}\left[\prod_{i=0}^{k}f_i(X_{t_i},V_{t_i})\right]\\
	=\langle h,T_{t_1,\infty}^{\mathbf{\Phi}}(f_1T_{t_2-t_1,\infty}^{\mathbf{\Phi}}(f_2...T_{t_{k-1}-t_{k-2},\infty}^{\mathbf{\Phi}}(f_{k-1}T_{t_k-t_{k-1},\infty}^{\mathbf{\Phi}}f_k)...))\rangle_{\mu_{\mathbf{\Phi}}}.	
	\end{equation}
	In particular, $\mathbb P_{h\mu_{\mathbf{\Phi}}}  $ solves the martingale problem for the generator $ (L_{\mathbf{\Phi},2},D(L_{\mathbf{\Phi},2})) $ of $ \sg{\mathbf \Phi}{,2} $, i.e.,
	denote by $ (X_t, V_t)_{t\geq0} $ the coordinate process on $ C([0,\infty),\{\Phi_1,\Phi_2<\infty\}) $. Then for $f\in D(L_{\mathbf{\Phi},2})$ the process $(M^{[ f ]}_t)_{t\geq0}$ defined by 
	\begin{equation}\label{Martingale2}
	M^{[ f ]}_t:=f(X_t, V_t)-f(X_0, V_0)-\int_{[0,t]}L_{\mathbf{\Phi},2}f(X_s, V_s)\,ds,\quad t\geq0,
	\end{equation}
	is a martingale w.r.t. the filtration $ \left(\mathcal F_t\right)_{t\geq0} $, $ \mathcal F_t=\sigma\left((X_s, V_s)\mid 0\leq s\leq t\right) $, and $\mathbb P_{h\mu_{\mathbf{\Phi}}} $. Additionally, if $ f^2\in D(L_{\mathbf{\Phi},2}) $ and $ L_{\mathbf{\Phi},2}f\in L^{4}(\Rd{2d},\mu_{\mathbf{\Phi}}) $ then the process $(N^{[ f ]}_t)_{t\geq0}$ defined by 
	\begin{equation}\label{QVP}
	N^{[ f ]}_t:=\left(M^{[ f ]}_t\right)^2-\int_{[0,t]}L_{\mathbf{\Phi},2}(f^2)(X_s, V_s)-2(fL_{\mathbf{\Phi},2}f)(X_s,V_s)\,ds,\quad t\geq0,
	\end{equation}
	is also a martingale w.r.t. $\mathbb P_{h\mu_{\mathbf{\Phi}}}  $ and the filtration $ \left(\mathcal F_t\right)_{t\geq0} $.
\end{T}

\begin{Rem}\label{RemarkMP}
	\mbox{}\vspace{-\topskip}
	\begin{enumerate}
		\item[(i)] Recall the situation of Theorem \ref{MP}. For $ f\in D(L^{(2)}_\mathbf{\Phi}) $ and $ 0\leq t\leq T<\infty $ the random variables in (\ref{Martingale2})
		%$ \int_{[0,T]}L_\mathbf{\Phi}f(X_s,V_s)\,ds $ 
		are well-defined, i.e., $ \mathbb P_{h\mu_{\mathbf{\Phi}}} $-a.s.~independent of the $ \mu_{ \mathbf\Phi} $ representative of $ f $ and $ L_{\mathbf{\Phi},2}f $, see \cite{CG10}[Lemma 5.1] for details. In particular it holds 
		\begin{equation*}
		\left\lVert\int_{[0,T]} \abs{L_{\mathbf{\Phi},2}f}(X_s,V_s)\,ds\right\rVert_{L^2(\mathbb P_{h\mu_{\mathbf{\Phi}}})}\leq T\norm{h}_{L^2(E,\mu)}\norm{L_{\mathbf{\Phi},2}f}_{L^2(E,\mu)}.
		\end{equation*}
		Hence, $ \int_{[0,T]} \abs{L_{\mathbf{\Phi}}f}(X_s,V_s)\,ds$ is finite $ \mathbb P_{h\mu_{\mathbf{\Phi}}} $-a.s.. On the negligible event 
		$$\bigcup_{T\in\mathbb N}\left\{\int_{[0,T]}\abs{L_{\mathbf{\Phi},2}f}(X_s,V_s)\,ds=\infty\right\} $$ we modify $ \int_{[0,t]} L_{\mathbf{\Phi},2}f(X_s,V_s)\,ds $ to be zero for all $ t\geq0 $ to obtain a continuous version of the process $ \left(\int_{[0,t]} L_{\mathbf{\Phi},2}f(X_s,V_s)\,ds\right)_{t\geq0} $. Hence, in the following we may assume that for continuous $ f $ the process  $(M^{[ f ]}_t)_{t\geq0}$ has continuous paths.
		\item[(ii)] The results from the previous Theorem also hold for the formal adjoint $ \hat L_{\mathbf\Phi} $, i.e., for $ h $ as in Theorem \ref{MP} there exists a law $\hat{\mathbb P}_{h\mu_{\mathbf{\Phi}}} $ on $C([0,\infty),\{\Phi_1,\Phi_2<\infty\}) $ with initial distribution $ h\mu_{\mathbf\Phi} $ which is associated with $ \asg{{\mathbf \Phi}}{,1} $ in the sense of (\ref{associated}), see \cite[Remark 3.3.]{CG10}. We use this fact later in the proof of Theorem \ref{Theoremtightness}.		
	\end{enumerate}
\end{Rem} 
% !TeX spellcheck = en_US
\section[Limit operator and limit process]{Limit operator and limit process}\label{LimitOperator}
%\subsection[M-Dissipativity of the limit operator $ L $]{M-Dissipativity of the limit operator $ L_0 $} 
This section consists of a brief summary of the functional analytic objects related to the overdamped Langevin equation (\ref{Overdamped1}) and the construction of martingale solutions for its generator. Denote by $ (B_t)_{t\geq o} $ a Brownian motion and recall the overdamped equation (\ref{Overdamped1}) 
\begin{equation}\label{Overdamped2}
dX^0_t=-\nabla\Phi_1(X^0_t)dt+\sqrt{2}dB_t.
\end{equation}
 The generator of (\ref{Overdamped2}) is given through 
\begin{equation}\label{limitop1}
L_{\Phi_1}f=\Delta f-\nabla\Phi_1\cdot\nabla f,  \quad f\in C^\infty_c(\{\Phi_1<\infty\}).
\end{equation}
 Recall the measure $ \mu_{{\Phi}_1} $ on $ \left(\Rd{d},\mathcal B(\Rd{d})\right) $ according to Notation \ref{mu}. Assuming $ (\Phi_12)-(\Phi_14)^2 $ one can use Proposition \ref{H1infty} to check that 
  the operator $ \left(L_{\Phi_1},C^\infty_c(\{\Phi_1<\infty\})\right) $ is symmetric and negative definite on the Hilbert space $ \mathcal H_{\Phi_1}=L^2(\Rd{d},\mu_{\Phi_1}) $, hence, closable. In particular, one can prove as in Corollary \ref{invariant} $ \int_{\Rd{d}}L_{\Phi_1}fd\mu_{\Phi_1}=0 $ for all $ f\in C^\infty_c(\{\Phi_1<\infty\}) $.
We make additional assumptions on $ \Phi_1 $.
\pagebreak
\begin{C}\label{C5}
 	\mbox{}\vspace{-\topskip}
 	\begin{itemize}[align=parleft, labelsep=0.5cm]
 		\item[$(\Phi_1\thePHIeins)$]\tabto{1cm}
 			The operator $\left(L_{\Phi_1},C^\infty_c(\{\Phi_1<\infty\})\right) $ is closable and its closure is the \\
 			\tabto{1cm}generator of a strongly continuous contraction semigroup $ \sg{\Phi_1}{,2} $ on $ \mathcal H_{\Phi_1}$.
 		\stepcounter{PHIeins}
 		\item[$(\Phi_1\thePHIeins)$]\tabto{1cm}
 			$ \mu_{ \Phi_1} $ is a finite measure, i.e., $\mu_{ \Phi_1}(\Rd{d})=\int_{\Rd{d}}e^{-\Phi_1}dx<\infty$.
 		\stepcounter{PHIeins}
 	\end{itemize}	
\end{C}

\begin{Rem}\label{RemSingPot}	
%	\mbox{}\vspace{-\topskip}
%	\begin{enumerate}
%		\item 
		The assumption $ (\Phi_15) $ still allows singular potentials $ \Phi_1 $. A very detailed discussion, including handy sufficient conditions and examples can be found in \cite[Section 4.2,4.3]{CG10}.
%		\item The assumption $(\Phi_16)$ is a sufficient condition for the process in Theorem \ref{processlimit} to have infinite life-time.
%	\end{enumerate}
	
\end{Rem}
\begin{T}\label{processlimit}
	Assume $ (\Phi_12),(\Phi_14)^2,(\Phi_15),(\Phi_16) $. Then the bilinear form $ \left(\mathcal E_{\Phi_1} ,C^\infty_c(\{\Phi_1<\infty\})\right) $ is closable and its closure $ \left(\mathcal E_{\Phi_1} ,D(\mathcal E_{\Phi_1} )\right) $ is a symmetric, quasi-regular Dirichlet form. Hence, there exists a $ \mu_{\Phi_1}$-tight special standard process
	\begin{equation*}\mathbb M_{\Phi_1} =\left(\Omega,\mathcal F , (\mathcal F _t)_{t\geq 0},(X _t)_{t\geq 0},(\mathbb P _x)_{x\in\{\Phi_1<\infty\}^\Delta}\right)
	\end{equation*}
	which is properly associated with  
	$ (\mathcal E_{\Phi_1} ,D(\mathcal E_{\Phi_1} )) $ in the resolvent sense.
	For each probability distribution $ \nu $ on $ \{\Phi_1<\infty\} $ being absolutely continuous w.r.t. $ \mu_{ \Phi_1} $ define the law $ \mathbb P_\nu(\cdot)=\int_{\{\Phi_1<\infty\}}\mathbb P_x(\cdot)d\nu(x) $. Then $ \mathbb P_\nu $-a.s. the paths are continuous and have infinite life-time.
	\begin{proof}
		Under the assumptions $ (\Phi_12),(\Phi_14)^2$ one obtains 
		\begin{equation}\label{preform}
		\mathcal E_{\Phi_1}(f,g)=-(L_{\Phi_1}f,g)_{\mathcal H_{\Phi_1}}, \quad f,g\in C^\infty_c(\{\Phi_1<\infty\}).
		\end{equation}
		Hence, the form $ (\mathcal E_{\Phi_1}, C^\infty_c(\{\Phi_1<\infty\})) $ is closable by \cite[Proposition I.3.3.]{MaR92}.
		The quasi-regularity of $ \left(\mathcal E_{\Phi_1} ,D(\mathcal E_{\Phi_1} )\right) $ holds by assumption $ (\Phi_15) $ and  \cite[IV.4.a)]{MaR92}. The sub-Markovian property of $ \sg{\Phi_1}{,2} $ can be proven as in \ref{maintheorem}, i.e., one shows $ \int_{\Rd{d}}L_{\Phi_1}fd\mu_{\Phi_1}=0 $ for all $ f\in  C^\infty_c(\{\Phi_1<\infty\}) $. Hence, \cite[Theorem IV.3.5]{MaR92} provides the existence of $ \mathbb M_{\Phi_1} $.
		Denote by $ \sg{\Phi_1}{,1} $, $ \sg{\Phi_1}{,\infty} $ the semigroups on $ L^1(\Rd{d},\mu_{\Phi_1}) $ and $ L^\infty(\Rd{d},\mu_{\Phi_1}) $, respectively, induced by the symmetric sub-Markovian semigroups $ \sg{\Phi_1}{,2} $, see \cite[Lemma 1.3.11.]{C10}. Denote by $\left(L^{(1)}_{\Phi_1},D\left(L^{(1)}_{\Phi_1}\right)\right)$ the generator of $ \sg{\Phi_1}{,1} $. Using \cite[Lemma 1.3.11.(iii)]{C10} and assumption $ (\Phi_16) $ one easily proves $ \int_{\Rd{d}}L^{(1)}_{\Phi_1}fd\mu_{\Phi_1}=0 $ for all $ f\in D\left(L^{(1)}_{\Phi_1}\right) $. Hence, $ \mu_{\Phi_1} $ is an invariant measure for the semigroup $ \sg{\Phi_1}{,1} $. Consequently, the semigroup $ \sg{\Phi_1}{,\infty} $ is conservative, see also the construction of $ \sg{\Phi_1}{,\infty} $ in \cite[Lemma 1.3.11.]{C10}.
		The continuity statement follows immediately by \cite[Theorem V.1.11.]{MaR92}.
	\end{proof}	
\end{T}
\noindent We obtain the analogous statement as in Theorem \ref{MP}.
% one can, due to the continuous sample paths of the process $ \mathbb M_0  $, construct for each $ h\in \mathcal H_0$ such that $ h $ is a probability density w.r.t. $ \mu_{\Phi_1} $ a probability measure $ \mathbb P_{h\mu_{\Phi_1}} $ on $ \left(C([0,\infty),\Rd{d}),\mathcal B(C([0,\infty),\Rd{d}))\right) $ which is associated to the semigroup $ \sg{0} $.
% The proof of the following corollary can be found in \cite[Lemma 4.10]{VN}. The fact about the martingale problem holds again due to \cite[Lemma 5.1]{CG10}.
\begin{Co}\label{MPL}
	Let  $ h\in L^1(\Rd{d},\mu_{\Phi_1})\cap L^2(\Rd{d},\mu_{\Phi_1})$ be a probability density w.r.t. $ \mu_{\Phi_1} $.
	Then there exists a probability law $\mathbb P_{h\mu_{\Phi_1}} $
	on $C([0,\infty),\{\Phi_1<\infty\}) $ with initial distribution $ h\mu_{\Phi_1} $ which is associated with the sub-Markovian strongly continuous contraction semigroup $ \sg{\Phi_1}{,2}$ in the sense that for all
	$ f_1,...,f_k\in L^\infty(\Rd{d},\mu_{\Phi_1}) $ and $ 0\leq t_1<...< t_k $, $ k\in\mathbb N $, it holds
	\begin{equation}
	\mathbb{E}\left[\prod_{i=0}^{k}f_i(X_{t_i})\right]\\
	=\langle h,T_{t_1,\infty}^{\Phi_1}(f_1T_{t_2-t_1,\infty}^{\Phi_1}(f_2...T_{t_{k-1}-t_{k-2},\infty}^{\Phi_1}(f_{k-1}T_{t_k-t_{k-1},\infty}^{\Phi_1}f_k)...))\rangle_{\mu_{\Phi_1}},
	\end{equation}
	where $ \mathbb{E} $ denotes integration w.r.t. $ \mathbb P_{h\mu_{\Phi_1}} $.
	In particular, the measure $\mathbb P_{h\mu_{\Phi_1}}  $ solves the martingale problem for the generator $ \left(L_{\Phi_1},D(L_{\Phi_1})\right) $.
	%	, i.e. for all $ f_1,...,f_k\in L^\infty(\Rd{2d},\mu_0) $ and $ 0\leq t_1<...< t_k $, $ k\in\mathbb N $, it holds
	%	\begin{align*}
	%	&\mathbb{E}\left[\prod_{i=0}^{k}f_i(X_{t_i},V_{t_i})\right]\\
	%	&=\langle h,T_{t_1,\infty}(f_1T_{t_2-t_1,\infty}(f_2...T_{t_{k-1}-t_{k-2},\infty}(f_{k-1}T_{t_k-t_{k-1},\infty}f_k)...))\rangle_{\mu_0}.
	%	\end{align*}
\end{Co}
\begin{Rem}
	\mbox{}\vspace{-\topskip}
\begin{enumerate}		
		\item One can prove stronger statements concerning life-time and continuity of the process $ \mathbb M_{\Phi_1} $. Since we only work in the following with laws $\mathbb P_{h\mu_{\Phi_1}} $ as in Corollary \ref{MPL} we restrict ourselves to the weaker statements.	
		\item In \cite{KR05} and the references therein strong solutions even for time-dependent and singular drifts of (\ref{Overdamped1}) are constructed. Under additional mild regularity assumptions on $ \Phi_1 $ we can show similar as below that weak solution can be constructed from the measure $ \mathbb P_{h\mu_{\Phi_1}} $ by proving e.g. that the functions $ f(x)=x_i $, $ i=1,...,d $ are contained in the domain $ D(L_{\Phi_1}) $.
	\end{enumerate}

\end{Rem}

\section[Velocity scaling and semigroup convergence]{Velocity scaling and semigroup convergence}\label{CoFuAnOb}
This section consists of a semigroup convergence result. 
For $ \varepsilon>0 $ we define a scaled velocity potential
\begin{equation}\label{scaling}
\Phi_2^\varepsilon(\cdot)=\Phi_2\left(\frac{\cdot}{\varepsilon}\right) +\ln(\varepsilon^d).
\end{equation} 
The constant $ \ln(\varepsilon^d) $ doesn't affect the generator and is only a renormalization constant. The assumptions $ (\Phi_21)-({\Phi}_27) $ hold true for $ \Phi_2^\varepsilon$ since they hold true for $ \Phi_2$.
Similar as before we write $ \mathbf{\Phi}^\varepsilon=\left(\Phi_1,\Phi_2^\varepsilon\right) $. We denote by $ \mu_\varepsilon $ the measure $ \mu_{\mathbf\Phi^\varepsilon} $.
 Hence, Theorem \ref{mdissiL1} and Theorem \ref{MP}
apply also for the operator $\left(L^1_{\mathbf{\Phi}^\varepsilon},C^\infty_c(\{\Phi_1,\Phi_2^\varepsilon<\infty\})\right)$
defined on $ L^1(\Rd{2d},\mu_\varepsilon) $ and its closure is denoted by $\left(L^1_{\mathbf{\Phi}^\varepsilon,1},D(L^1_{\mathbf{\Phi}^\varepsilon,1})\right)$.
%which we abbreviate again by $ L_{\mathbf{\Phi}^\varepsilon} $.
Furthermore, we obtain a strongly continuous contraction semigroups $ \sg{\varepsilon}{,2}=\sg{\mathbf{\Phi}^\varepsilon}{,2} $ on the Hilbert space $ \mathcal H_\varepsilon= L^2(\Rd{2d},\mu_\varepsilon)$, see Lemma \ref{Lpgen} and its previous discussion. The generator $\left(L^1_{\mathbf{\Phi}^\varepsilon,2},D(L^1_{\mathbf{\Phi}^\varepsilon,2})\right)$ of $ \sg{\varepsilon}{,2} $ we abbreviate by $ (L_\varepsilon,D(L_\varepsilon)) $. Observe that $ (L_\varepsilon,D(L_\varepsilon)) $ is an extension of $\left(L^1_{\mathbf{\Phi}^\varepsilon},C^\infty_c(\{\Phi_1,\Phi_2^\varepsilon<\infty\})\right)$ considered as an operator on $ \mathcal H_\varepsilon$. Additionally we denote by $ \sg{0}{} $ the semigroup $ \sg{\Phi_1}{} $ on $ \mathcal H_0:=\mathcal H_{\Phi_1} $.
In the following we show convergence of the Hilbert spaces $ \mathcal H_\varepsilon $ towards the Hilbert space $ \mathcal H_0 $ from Section \ref{LimitOperator} in the sense of Kuwae-Shioya, i.e., there exists a dense subset $ \mathcal C $ of $ \mathcal H_0 $ and for every $ \varepsilon>0 $ there exists a linear map 
\begin{equation}\label{Psieps1}
\Psi_\varepsilon:\mathcal C\longrightarrow \mathcal H_\varepsilon,
\end{equation}
such that 
\begin{equation}\label{HSconvergence}
\lim_{\varepsilon\to0}\norm{\Psi_\varepsilon(u)}_{\mathcal H_\varepsilon}=\norm{u}_{\mathcal H_0},\text{ for all } u\in\mathcal C.
\end{equation}
If (\ref{HSconvergence}) holds we say that the family of Hilbert spaces $ \left(\mathcal H_\varepsilon\right)_{\varepsilon>0} $ converges to $ \mathcal H_0 $ along the family $ \left(\Psi_\varepsilon\right)_{\varepsilon>0} $ and we use the short hand notation $\mathcal{H}_\varepsilon\xrightarrow{\left(\Psi_\varepsilon\right)_{\varepsilon>0}}\mathcal{H}_0$. In this case we say that $ f_\varepsilon\in\mathcal H_\varepsilon $, $ \varepsilon>0 $, converges to $ f\in\mathcal H_0 $ (Notation: $ f_\varepsilon\longrightarrow f $ along $\mathcal{H}_\varepsilon\xrightarrow{\left(\Psi_\varepsilon\right)_{\varepsilon>0}}\mathcal{H}_0$) if 
\begin{align}
\norm{f_\varepsilon}_{\mathcal H_\varepsilon}&\xrightarrow{\varepsilon\to0}\norm{f}_{\mathcal H_0}\label{normconv}\\
 (f_\varepsilon,\Psi_\varepsilon(\varphi))_{\mathcal H_\varepsilon}&\xrightarrow{\varepsilon\to0}(f,\varphi)_{\mathcal H_0}  \text{ for all } \varphi\in \mathcal C\label{weakconv}. 
\end{align}

Furthermore, we prove convergence of the semigroups $ \sg{\varepsilon}{,2} $, $ \varepsilon>0 $, towards the semigroup $ \sg{0}{} $ along $\mathcal{H}_\varepsilon\xrightarrow{\left(\Psi_\varepsilon\right)_{\varepsilon>0}}\mathcal{H}_0$, i.e., for all $ t\geq0 $ it holds
\begin{equation}\label{sgc}
 f_\varepsilon\longrightarrow f  \text{ along } \mathcal{H}_\varepsilon\xrightarrow{\left(\Psi_\varepsilon\right)_{\varepsilon>0}}\mathcal{H}_0 \text{ implies } T^\varepsilon_{t,2}f_\varepsilon\longrightarrow T^0_{t,2}f  \text{ along } \mathcal{H}_\varepsilon\xrightarrow{\left(\Psi_\varepsilon\right)_{\varepsilon>0}}\mathcal{H}_0.
\end{equation}
To this end, we assume that $ \Phi_1$ and $\Phi_2 $, respectively, fulfill the additional assumptions:
\begin{C}\label{C6}
	\mbox{}\vspace{-\topskip}
	\begin{itemize}[align=parleft, labelsep=0.5cm] 	
		\item[$(\Phi_1\thePHIeins)$]\label{phi12}
		\stepcounter{PHIeins}\tabto{1cm}
		The measurable function $\abs{\nabla\Phi_1}$ is square integrable w.r.t. $\mu_{{\Phi}_1}$, i.e.,\tabto{1cm}
		{$
		\int_{\Rd{d}}\abs{\nabla\Phi_1}^2\,d\mu_{{\Phi}_1}=
		\int_{\Rd{d}}\abs{\nabla\Phi_1}^2e^{-\Phi_1}\,dx<\infty.
		$}
	\end{itemize}
\end{C}
\begin{C}\label{C7}
	\mbox{}\vspace{-\topskip}
	\begin{itemize}[align=parleft, labelsep=0.5cm] 
		\item[$(\Phi_2\thePHIzwei)$]\tabto{1cm}
		$ \Phi_2 $ has no singlarities, i.e., $ \{\Phi_2=\infty\}=\emptyset $.
		\stepcounter{PHIzwei} 
	\end{itemize}
\end{C}
\noindent Due to $(\Phi_27)$ we can assume $ \mu_{{\Phi}_2}(\Rd{d})=1 $. Furthermore, we define the following maps $p_x,p_v,\sigma:\Rd{2d}\longrightarrow\Rd{d}$, where $ \sigma(x,v)=x+v$, $p_x(x,v)=x$, $p_v(x,v)=v $.
Next we define the maps $ \Psi_\varepsilon $ from (\ref{Psieps1}).
\begin{D}\label{Psieps2}
	Let $\varepsilon>0$ and choose a symmetric cut off function $\eta_\varepsilon\in C^\infty_c(\Rd{d})$, s.t. 
	\begin{enumerate}\label{etaprop}
		\item[(i)] $\eta_\varepsilon(v)=\eta_\varepsilon(-v)$, for all $v\in \Rd{d}$, 
		$\eta_\varepsilon\equiv1$ on $B_{\varepsilon^{-2}}(0)$ and
		$\text{supp}(\eta_\varepsilon)\subseteq B_{2\varepsilon^{-2}}(0)$,
		\item[(ii)] $\abs{\nabla\eta_\varepsilon}\leq C\varepsilon^2$ and $\abs{\Delta\eta_\varepsilon}\leq C
		\varepsilon^4$, for a finite constant $C$ independent of $\varepsilon$.
	\end{enumerate}
	 We choose $ \mathcal{C}=C^\infty_c(\{\Phi_1<\infty\}) $ and define the convergence determining function $\Psi_\varepsilon$ by 
	\begin{equation}\label{Psi}
	\Psi_\varepsilon:\mathcal{C}\longrightarrow\mathcal{H}_\varepsilon,f\mapsto (f\circ\sigma)(\eta_{\varepsilon}\circ p_v).   
	\end{equation}
	Due to Proposition \ref{extensionweak} and Lemma \ref{Lpgen}(i) it holds $ \Psi_\varepsilon(\mathcal C)\subseteq C^\infty_c(\Rd{2d})\subseteq D(L_\varepsilon) $.
\end{D}
\begin{T}\label{CoSe}
	Assume $ (\Phi_12)-(\Phi_14)^2,(\Phi_15),(\Phi_17) $, $ (\Phi_21)-(\Phi_29) $ and one of the additional assumptions $ (i) $, $ (ii) $ of Theorem \ref{mdissiL1} to hold true. Then it holds, the family of Hilbert spaces $ \left(\mathcal H_\varepsilon\right)_{\varepsilon>0} $ converges along the family $\left(\Psi_\varepsilon\right)_{\varepsilon>0} $ defined in (\ref{Psi}) towards the Hilbert space $ \mathcal H_0 $ as $ \varepsilon $ tends to zero in the Kuwae-Shioya sense. Furthermore, the semigroups $ \sg{\varepsilon}{,2} $,  $\varepsilon>0$, converge towards $ \sg{0}{} $ along $\mathcal{H}_\varepsilon\xrightarrow{\left(\Psi_\varepsilon\right)_{\varepsilon>0}}\mathcal{H}_0$, i.e., (\ref{sgc}) holds.
	\begin{proof}
		We proceed as in \cite[Proposition 3.21., Theorem 3.22.]{VN}, where the special case $ \Phi_2(v)=\frac{1}{2}\abs{v}^2 $ is considered. For sake of completeness we give a short proof.
		For $ f \in \mathcal C $ we have to show $\norm{\Psi_\varepsilon f}_{\mathcal H_\varepsilon}\xrightarrow{\varepsilon\to0}\norm{f}_{\mathcal H_0}$.
		Using the symmetry of $ \eta_{\varepsilon} $ and $ \Phi_2 $ together with the transformation $ (x,v)\mapsto (x,-v) $ we rewrite the norm using the convolution $ \ast $, i.e., 
		\begin{equation}\label{normconv1}
		\norm{\Psi_\varepsilon f}_{\varepsilon}^2=	\int_{\mathbb{R}^{d}}f^2\ast(\eta_\varepsilon^2e^{-\Phi_2^\varepsilon})(x)e^{-\Phi_1}(x)\,dx.
		\end{equation}
		%As in \cite[Lemma 3.19]{VN} one defines 
		For $\alpha_\varepsilon:=\int_{\mathbb{R}^d}\eta_\varepsilon^2e^{-\Phi_2^\varepsilon}(v)\,dv$
		one can show $\alpha_\varepsilon\overset{\varepsilon\to 0}{\longrightarrow} 1 $, hence  $\left(\alpha_\varepsilon^{-1}\eta_\varepsilon^2e^{-\Phi_2^\varepsilon}\right)_{\varepsilon>0}$ is an approximate identity. Since $ f^2\in L^1(\Rd{d}) $ and $ e^{-\Phi_1}\in L^\infty(\Rd{d}) $ due to assumption $ ({\Phi}_12) $ the H"{o}lder inequality implies the desired result.
		
		Next we prove convergence of the semigroups generated by 
		$\left(L_\varepsilon,D(L_\varepsilon)\right)$ in $\mathcal{H}_\varepsilon$. Recall that the limit semigroup $\sg{0}{}$ has the closure of $ (L_{\Phi_1},C^\infty_c(\{\Phi_1<\infty\})) $ as its generator.
		We use that semigroup convergence is equivalent to convergence of the generators and in particular it suffices to have convergence of the generators on a core for the limit generator, i.e., we use \cite[Theorem 1.5.13]{C10}, \cite[Corollary 1.5.14]{C10}. 
%		In particular, the set $\mathcal{C}$ is a core for generator of $\sg{0}{}$.
		 Hence for $f\in \mathcal{C}=C^\infty_c(\{\Phi_1<\infty\})$ it suffices to show $(L_\varepsilon\Psi_\varepsilon f)_{\varepsilon>0}\longrightarrow L_0 f$ along 
		$\mathcal{H}_\varepsilon\xrightarrow{\left(\Psi_\varepsilon\right)_{\varepsilon>0}}\mathcal{H}_0$. 
		Let $ f:\Rd{d}\longrightarrow \Rd{}$ be smooth and $ i\in\{1,...,d\} $. Observe that the function $f\circ \sigma$  fulfills
		$\partial_{x_i} (f\circ \sigma)=\partial_i f\circ \sigma=\partial_{v_i} (f\circ \sigma)$. 
		We start with computing the expression $L_\varepsilon\Psi_\varepsilon f$ explicitly.
		According the previous observation we obtain
		\begin{align}
			L_\varepsilon\Psi_\varepsilon f=&(\Delta f\circ \sigma)( \eta_\varepsilon\circ p_v) + (f\circ \sigma)(\Delta \eta_\varepsilon\circ p_v)
			+2(\nabla f\circ\sigma)\cdot (\nabla\eta_\varepsilon\circ p_v)\notag\\
			&-\left(\nabla_v\Phi_2^\varepsilon\cdot(\nabla\eta_\varepsilon\circ p_v)\right) (f\circ\sigma)
			-(\nabla_x\Phi_1\cdot(\nabla f\circ\sigma)) \eta_\varepsilon\circ p_v\notag\\
			&-\left(\nabla_x\Phi_1\cdot(\nabla\eta_\varepsilon\circ p_v) \right) (f\circ\sigma).\label{Leps}
			\intertext{The aim is to establish that (\ref{Leps}) converges along 
				$\mathcal{H}_\varepsilon\xrightarrow{\left(\Psi_\varepsilon\right)_{\varepsilon>0}}\mathcal{H}_0$ towards}
			L_0 f=&\Delta f-\nabla\Phi_1\cdot\nabla f.
		\end{align}
		Since convergence along $\mathcal{H}_\varepsilon\xrightarrow{\left(\Psi_\varepsilon\right)_{\varepsilon>0}}\mathcal{H}_0$ is linear
		(see \cite[Lemma 2.1. (3)]{KS03}) it suffices to show convergence of the single summands in (\ref{Leps}), i.e., one shows
		
		\begin{tabular}{lclccc}
			1.&$(f\circ \sigma) (\Delta \eta_\varepsilon\circ p_v)$&\rdelim\}{4}{0pt}&\multirow{4}{*}{$\longrightarrow$}&\multirow{4}{*}{0}& \multirow{6}{*}{$\text{along }\mathcal{H}_\varepsilon\xrightarrow{\left(\Psi_\varepsilon\right)_{\varepsilon>0}}\mathcal{H}_0.$}\\
			2.&$(\nabla f\circ\sigma)\cdot (\nabla\eta_\varepsilon\circ p_v)$& & & &\\
			3.&$\left(\nabla_x\Phi_1\cdot(\nabla\eta_\varepsilon\circ p_v) \right) (f\circ\sigma)$& & & &\\
			4.&$\left(\nabla_v\Phi_2^\varepsilon\cdot(\nabla\eta_\varepsilon\circ p_v)\right) (f\circ\sigma)$& & & &\\
			5.&$ (\Delta f\circ \sigma)( \eta_\varepsilon\circ p_v) $& & $ \longrightarrow $ & $ \Delta f$  &\\
			6.&$ (\nabla_x\Phi_1\cdot(\nabla f\circ\sigma)) (\eta_\varepsilon\circ p_v) $ & & $  \longrightarrow $ & $\nabla\Phi_1\cdot\nabla f $&
		\end{tabular}
		\vspace{0.5cm}
		\newline
		To prove convergence in 1.-4. one 
		%uses \cite[Lemma 2.1. (1)]{KS03} and 
		checks that the respective norms of the elements converge to zero, see \cite[Lemma 2.1. (1)]{KS03}. But this holds due to the choice of $ \eta_\varepsilon $ and a convolution argument as in (\ref{normconv1}).
		The statements in 5. and 6. are obtained by the same convolution argument.
		Taking 1.-6. together we obtain 
		\begin{equation}
		L_\varepsilon\Psi_\varepsilon f\longrightarrow L_0f\text{ along }
		\mathcal{H}_\varepsilon\xrightarrow{\left(\Psi_\varepsilon\right)_
			{\varepsilon>0}}\mathcal{H}_0,\quad\forall f\in C^\infty_c(\{\Phi_1<\infty\})
		\end{equation}
	\end{proof}	
\end{T}

% !TeX spellcheck = en_US

\section[Convergence in law of weak solutions]{Convergence in law of weak solutions}\label{WCoL}
 Throughout this section let $ \varepsilon>0 $ and $ h_\varepsilon\in\mathcal H_\varepsilon $ and $ h_0\in\mathcal H_0 $ be probability densities w.r.t.
 $ \mu_\varepsilon $ and $ \mu_0:=\mu_{ \Phi_1} $, respectively. Furthermore, let $ \mathbb P_{h_\varepsilon\mu_\varepsilon} $ by the martingale solution for $ \left(L^1_{\mathbf\Phi^\varepsilon,2},D\left(L^1_{\mathbf\Phi^\varepsilon,2}\right)\right) $ with initial distribution $ h_\varepsilon\mu_\varepsilon $ given by Theorem \ref{MP} and $ \mathbb P_{h_0\mu_0} $ be the measure from Corollary \ref{MPL}. The measures $ \mathbb P_{h_\varepsilon\mu_\varepsilon} $ and $ \mathbb P_{h_0\mu_0} $ are defined on $ C\left([0,\infty),\{\Phi_1<\infty\}\times\Rd{d}\right) $ and $ C\left([0,\infty),\{\Phi_1<\infty\}\right) $, respectively.
 In the following we consider them as measures on $ C\left([0,\infty),\Rd{2d}\right) $ and $ C\left([0,\infty),\Rd{d}\right) $. Indeed, we consider the continuous embeddings
 \begin{align*}
 i_{2d}&:C\left([0,\infty),\{\Phi_1<\infty\}\times\Rd{d}\right)\longrightarrow C\left([0,\infty),\Rd{2d}\right),\omega\mapsto\omega,\\
 i_d&:C\left([0,\infty),\{\Phi_1<\infty\}\right)\longrightarrow C\left([0,\infty),\Rd{d}\right),\omega\mapsto\omega.
 \end{align*}
 We also denote by  $ \mathbb P_{h_\varepsilon\mu_\varepsilon} $ and $ \mathbb P_{h_0\mu_0} $ the pushforwards $ \mathbb P_{h_\varepsilon\mu_\varepsilon}\circ i_{2d}^{-1} $ and $ \mathbb P_{h_0\mu_0}\circ i_d^{-1} $, respectively, to ease the notation. Observe that these measures are still associated with the respective semigroup.
Additionally, we define the continuous coordinate projection
 \begin{equation}\label{ProjX}
\noindent P_X:C\left([0,\infty),\Rd{2d}\right)\longrightarrow C\left([0,\infty),\Rd{d}\right),(x_t,v_t)_{t\geq 0} \mapsto (x_t)_{t\geq0}.
 \end{equation}
 In this section we prove weak convergence of $\mathbb P_{h_\varepsilon\mu_\varepsilon}^X:= \mathbb P_{h_\varepsilon\mu_\varepsilon}\circ P_X^{-1}$ towards $ \mathbb P_{h_0\mu_0} $ for $ \varepsilon\to 0 $ as measures on $ C\left([0,\infty),\Rd{d}\right) $.
 At first, weak convergence of the finite dimensional distributions (f.d.d.) is shown via the convergence of the associated semigroups $ \sg{\varepsilon}{,2} $, i.e., Theorem \ref{CoSe}. In a second step we prove tightness implying weak convergence.
\begin{T}\label{fdd}
Assume $ (\Phi_12)-(\Phi_14)^2,(\Phi_15)-(\Phi_17)$ and $ (\Phi_21)-(\Phi_29)$.
%Let $ h_\varepsilon\in \mathcal H_\varepsilon$, $ \varepsilon>0 $, $ h_0 \in \mathcal H_0 $ be probability densities w.r.t. $ \mu_\varepsilon $, $ \varepsilon>0 $, and $ \mu_0 $, respectively.
 If $h_\varepsilon\mu_\varepsilon $ converges weakly to $ h_0\mu_0\otimes \delta_0 $, where $ \delta_0 $ is the Dirac measure in zero on $ \Rd{d} $, as measures on $ \Rd{2d} $ and $ \sup_{\varepsilon>0}\norm{h_\varepsilon}_{L^2(\mu_\varepsilon)}<\infty $ then the  f.d.d. of $ \mathbb P_{h_\varepsilon\mu_\varepsilon}^X $ converge weakly to the f.d.d. of $ \mathbb P_{h_0\mu_0} $ as $ \varepsilon\to0 $.
	\begin{proof}
		Let $ \left(X_t\right)_{t\geq0} $ and $ \left(X_t,V_t\right)_{t\geq0} $ be the coordinate processes on $ C\left([0,\infty),\Rd{d}\right) $ and $ C\left([0,\infty),\Rd{2d}\right)$, respectively. Then it holds $ X_t\circ P_X=p_x\circ\left(X_t,V_t\right) $ for all $ t\geq0$.
		Let $ 0\leq t_1<...< t_k  $, $ k\in\mathbb N $ and define 
		$ \mathbb P_{h_\varepsilon\mu_\varepsilon}^{X,t_1,...,t_k}:= \mathbb P_{h_\varepsilon\mu_\varepsilon}^X\circ\left(X_{t_1},...,X_{t_k}\right)^{-1}$ and  $\mathbb P_{h_0\mu_0}^{t_1,...,t_k}:=\mathbb P_{h_0\mu_0}\circ \left(X_{t_1},...,X_{t_k}\right)^{-1}$. Additionally, let $ F:\Rd{dk}\longrightarrow\Rd{}$ be of the form $ F(x_1,...,x_k)=\prod_{i=1}^kf_i(x_i) $, $ f_i\in C^\infty_c(\Rd{d}) $, $ i=1,...,k $.
		 By the association of $ \mathbb P_{h_\varepsilon\mu_\varepsilon} $ with $ \sg{\mathbf\Phi^\varepsilon}{,1} $ and $T^\varepsilon_{t,2}=T^{\mathbf\Phi^\varepsilon}_{t,\infty}$ on $ L^2(\Rd{2d},\mu_\varepsilon)\cap L^\infty(\Rd{2d},\mu_\varepsilon) $ it holds
		\begin{equation}
		\int_{\Rd{dk}}F\,d\mathbb{P}_{h_\varepsilon\mu_\varepsilon}^{X,t_1,...,t_k}
		=\int_{\Rd{2d}}h_\varepsilon\underbrace{T_{t_1,2}^\varepsilon(f_1\circ p_xT_{t_2-t_1,2}^\varepsilon (f_2\circ p_x...T_{t_k-t_{k-1},2}^\varepsilon f_k\circ p_x))...)}_{F_\varepsilon^{t_1,...,t_k}}\,d\mu_\varepsilon.
		\end{equation}
		Observe that for $ g\in C^\infty_c(\Rd{d}) $ the constant sequence $ g\circ p_x\in \mathcal H_\varepsilon $ converges to $ g $ along $\mathcal{H}_\varepsilon\xrightarrow{\left(\Psi_\varepsilon\right)_{\varepsilon>0}}\mathcal{H}_0$. Furthermore, for $ f_\varepsilon\longrightarrow f $ along $\mathcal{H}_\varepsilon\xrightarrow{\left(\Psi_\varepsilon\right)_{\varepsilon>0}}\mathcal{H}_0$ it holds $ (g\circ p_x) f_\varepsilon\longrightarrow g f $.
		Applying Theorem \ref{CoSe} and the previous observations inductively we see that $ F_\varepsilon^{t_1,...,t_k} $ converges to $F_0^{t_1,...,t_k}:= T^0_{t_1}(f_1T^0_{t_2-t_1}(f_2...T^0_{t_k-t_{k-1}}f_k)...) $ along $ \mathcal{H}_\varepsilon\xrightarrow{\left(\Psi_\varepsilon\right)_{\varepsilon>0}}\mathcal{H}_0 $. Furthermore, the densities $h_\varepsilon $ converge weakly towards $ h_0 $ along $ \mathcal{H}_\varepsilon\xrightarrow{\left(\Psi_\varepsilon\right)_{\varepsilon>0}}\mathcal{H}_0 $ by \cite{T06}[Lemma 2.13]. We conclude
		\begin{equation*}
		\int_{\Rd{dk}}F\,d\mathbb{P}_{h_\varepsilon\mu_\varepsilon}^{X,t_1,...,t_k}=\left(h_\varepsilon,F_\varepsilon^{t_1,...,t_k}\right)_\varepsilon\xrightarrow{\varepsilon\to0}\left(h_0,F_0^{t_1,...,t_k}\right)_0=\int_{\Rd{dk}}F\,d\mathbb{P}_{h_0\mu_0}^{t_1,...,t_k}.
		\end{equation*}		
		Since the functions $ F $ of this kind are strongly separating
		\cite[Chapter 3, Theorem 4.5]{EK86} yields the claim.
	\end{proof}
\end{T}
\noindent To prove tightness we choose an appropriate metric $ m $ on our state space
$ \Rd{2d} $ inducing the euclidean topology. Let $ i\in\{1,..,d\} $ and define the functions $ f_i$, $g_i $ in the following way:
\begin{equation}\label{fg}
f_i:\Rd{2d}\longrightarrow \Rd{},(x,v)\mapsto x_i+v_i,\quad
g_i:\Rd{2d}\longrightarrow \Rd{},(x,v)\mapsto v_i.
\end{equation}
Let the metric $ m $ on $\Rd{2d}$ be given by
\begin{equation}\label{metric}
m((x,v),(\tilde x,\tilde v))=\sum_{i=1}^d\abs{f_i((x,v))-f_i((\tilde x,\tilde v))}+\abs{g_i((x,v))-g_i((\tilde x,\tilde v))}.
\end{equation}
We need further assumptions on 
%assume Assumption \ref{C8} and \ref{C9} to hold true for 
$\Phi_1 $ and $ \Phi_2 $, respectively.
\pagebreak
{\begin{C}\label{C8}
	\mbox{}\vspace{-\topskip}
	\begin{itemize}[align=parleft, labelsep=0.5cm]
		\item[$(\Phi_1\thePHIeins)$] \tabto{1cm}
		$ \int_{\mathbb{R}^{d}}\abs{x}^{2k}e^{-\Phi_1}\,dx<\infty, \quad k=1,2$
		\stepcounter{PHIeins}
		\item[$(\Phi_1\thePHIeins)$] \tabto{1cm}
		$ \int_{\mathbb{R}^{d}}\abs{\nabla\Phi_1}^{4}e^{-\Phi_1}\,dx<\infty.$
		\stepcounter{PHIeins}
	\end{itemize}	
\end{C}}
\begin{C}\label{C9}
	\mbox{}\vspace{-\topskip}
	\begin{itemize}[align=parleft, labelsep=0.5cm]
		\item[$(\Phi_2\thePHIzwei)$] \tabto{1cm}
		$ \int_{\mathbb{R}^{d}}\abs{v}^{2k}e^{-\Phi_2}\,dv<\infty, \quad k=1,2$
		\stepcounter{PHIzwei}
		\item[$(\Phi_2\thePHIzwei)$]  \tabto{1cm}
		$ \int_{\mathbb{R}^{d}}\abs{\nabla\Phi_2}^{4}e^{-\Phi_2}\,dv<\infty.$
		\stepcounter{PHIzwei}
	\end{itemize}	
\end{C}
\noindent Due to $ (\Phi_16) $ and $ (\Phi_27) $ the measure $ \mu_{\mathbf{\Phi}} $ is finite, hence, w.l.o.g. we assume that $ \mu_{\varepsilon} $ is a probability measure for all $ \varepsilon $.
% To prove tightness we first establish tightness for the family $ \left(\mathbb P_{\mu_\varepsilon}\right)_{\varepsilon>0} $, i.e. we first consider the case $ h_\varepsilon=1 $ for all $ \varepsilon>0 $. 
%  is tight implying tightness of $ \left(\mathbb P^X_{\mu_\varepsilon}\right)_{\varepsilon>0} $.
For $ h_\varepsilon=1 $ the measure $ \mu_\varepsilon $ is invariant for $ \mathbb P_{\mu_\varepsilon} $ for all $ \varepsilon>0 $, i.e., the one dimensional distributions of $ \mathbb P_{\mu_\varepsilon} $ are given by $ \mu_\varepsilon $. Furthermore, the family $ \mu_\varepsilon $, $ 0<\varepsilon\leq 1 $, is tight.
 Denote by $ (\hat L_\varepsilon,D(\hat L_\varepsilon)) $ the generator of the adjoint semigroup $ \asg{\varepsilon}{,2} $. 
\begin{Le}\label{figi}
	Assume $ (\Phi_12),(\Phi_13),(\Phi_15)-(\Phi_19) $ and $ (\Phi_21)-(\Phi_27),(\Phi_29)-(\Phi_211) $.
	For the functions $ f_i,g_i $, $ i\in\{1,..,d\} $, defined in (\ref{fg}) it holds $f_i,f_i^2,g_i,g_i^2\in  D(L_\varepsilon)\cap D(\hat L_\varepsilon) $ and 	
	\begin{align}
	L_\varepsilon f_i&= -\partial_{x_i}\Phi_1, & L_\varepsilon f_i^2&=2+2f_iL_\varepsilon f_i \label{fi}\\
	L_\varepsilon g_i&=-\partial_{v_i}\Phi_2^\varepsilon-\partial_{x_i}\Phi_1, &L_\varepsilon g_i^2&=2+2g_iL_\varepsilon g_i,\label{gi1}\\	
	\hat L_\varepsilon g_i&=-\partial_{v_i}\Phi_2^\varepsilon+\partial_{x_i}\Phi_1, &\hat L_\varepsilon g_i^2&=2+2g_i\hat L_\varepsilon g_i.\label{gi2}		  
	\end{align}
	\begin{proof}
		Due to Proposition \ref{extensionweak} and Lemma \ref{Lpgen}(i) we know that $ C^\infty_c(\Rd{2d}) $ is contained in $ D(L_\varepsilon)\cap D(\hat L_\varepsilon) $. 
		The assertions follow using suitable cut off functions. 
	\end{proof}
%	\begin{proof}
%		We only indicate how to prove the two statements from (\ref{gi}). 
%		Choose a family of cut off functions $ \left(\eta_n\right)_{n\in\mathbb N}\subseteq C^\infty_c(\Rd{d}) $ such that $\eta_n\equiv1$ on $B_{n}(0)$ and
%		$\text{supp}(\eta_n)\subseteq B_{2n}(0)$, $\abs{\partial^\alpha\eta_n}\leq \frac{C}{n^{\abs{\alpha}}}$ for all $ \alpha\in\mathbb N_0^d $, $ \abs{\alpha}\leq 2 $ and a positive constant $C$ independent of $n$. Then one defines $ g_i^n:=g_i(\eta_n\otimes\eta_n)\in C^\infty_c(\Rd{2d}) $. It holds $ g_i^n\longrightarrow g_i $, $ n\to\infty $ in $ L^2(\mu_{\varepsilon}) $ and in $ L^4(\mu_{\varepsilon}) $.
%%		 Since $ g_i^n$ is an element of $C^\infty_c(\Rd{2d}) $ we have an explicit representation of $ L_\varepsilon g_i^n $.
%		Explicit computation leads to $L_\varepsilon g_i^n\longrightarrow -\frac{1}{\sqrt{\varepsilon}}\partial_{v_i}\Phi_2\left(\frac{\cdot}{\sqrt{\varepsilon}}\right)-\partial_{x_i}\Phi_1 $, as $ n \to\infty $  in $ L^2(\mu_{\varepsilon}) $ and in $ L^4(\mu_{\varepsilon}) $, which proves the first assertion. For the second one we approximate $ g_i^2 $ by $ (g_i^n)^2 $. Then $ L_\varepsilon (g_i^n)^2 = 2 g_i^nL_\varepsilon g_i^n+ 2\abs{\nabla_vg_i^n}^2\longrightarrow  2 g_iL_\varepsilon g_i+ 2 $ in $ L^2(\mu_{\varepsilon}) $ as  $ n\to\infty $.
%	\end{proof}
\end{Le}
\begin{Rem}\label{weaksolution}
	Observe that the assumptions of the previous lemma imply that
	the coordinate process $ \left(X_t,V_t\right)_{t\geq0} $ on $ C\left([0,\infty),\Rd{2d}\right) $ is a weak solution  
%	 the measures $ \mathbb P_{h_\varepsilon\mu_{\varepsilon}} $ constitute weak solutions
	  to (\ref{sHseps1}), (\ref{sHseps2}) for $ \Phi_2^\varepsilon $ instead of $ \Phi_2 $ and $ \varepsilon=1 $ with initial distributions $ h_\varepsilon\mu_{\varepsilon} $ under $ \mathbb P_{h_\varepsilon\mu_{\varepsilon}} $. Indeed, let $ i\in\{1,...,d\}. $ Due to Lemma \ref{figi} we know that the function $g_i $ is in $ D(L_\varepsilon) $. By (\ref{QVP}) we know that the quadratic cross-variations of the continuous $ d- $dimensional martingale $ \left(M^{[g_i],\varepsilon}_t\right)^{i=1,..,d}_{t\geq0} $ is given by \begin{equation*}
	  \left<M^{[g_i],\varepsilon},M^{[g_j],\varepsilon}\right>_t=\delta_{ij}t,
	  \end{equation*}
	  where $ \delta_{ij} $ denotes the Kronecker delta.
	   Using L\'{e}vy's characterization of Brownian motion, we see that $ \left(M^{[g_i],\varepsilon}_t\right)^{i=1,..,d}_{t\geq0} $ is $ \sqrt{2} $ times a $ d- $dimensional Brownian motion. Computing the quadratic variation of $ \left(M^{[f_i-g_i],\varepsilon}_t\right)^{i=1,..,d}_{t\geq0} $ we obtain $ M^{[f_i-g_i],\varepsilon}_t=0 $ for all $ t\geq0 $. Hence, by comparing (\ref{sHseps1}), (\ref{sHseps2}) with (\ref{Martingale2}) for $ f_i-g_i $ and $ g_i $ we constructed a $ d- $dimensional Brownian motion $ (B_t)_{t\geq0} $ and a stochastic process $ (X_t,V_t)_{t\geq0} $ such that (\ref{sHseps1}), (\ref{sHseps2}) holds.  
\end{Rem}
%\begin{Le}\label{timerev}
%	Let $ \varepsilon>0 $  and $ T>0 $. The time reversed law $ \mathbb P^{T}_{\mu_\varepsilon}\circ r_T^{-1} $ is associated with $ \left(\hat{T}^\varepsilon_t\right)_{t\in[0,T]}$ in the sense of Definition \ref{SAs}. In particular, it holds
%	\begin{equation}
%	\hat{\mathbb P}^{T}_{\mu_\varepsilon}=\mathbb P^{T}_{\mu_\varepsilon}\circ r_T^{-1}.
%	\end{equation}
%	\begin{proof}
%		See \cite[Lemma 3.9]{GS16}.
%	\end{proof}
%\end{Le}
%\noindent Now we prove tightness:
%Unfortunately the space $C\left([0,\infty),\{\Phi_1<\infty\}\times\Rd{d}\right) $ is not appropriate to prove tightness on
%\noindent	Since $ \{\Phi_1<\infty\}\times\Rd{d} $ is merely open in $ \Rd{2d} $, we don't have a characterization of compactness in $C\left([0,\infty),\{\Phi_1<\infty\}\times\Rd{d}\right) $ in terms of the Arzela-Ascoli theorem. Hence, we will consider the family of image measures of $\left(\mathbb{P}_{\mu_\varepsilon}\right)_{\varepsilon>0} $ under the embedding
%\begin{equation*}\label{embedding}
%i:C\left([0,\infty),\{\Phi_1<\infty\}\times\Rd{d}\right)\longrightarrow C\left([0,\infty),\Rd{2d}\right),
%\end{equation*}
%which we also denote by $\left(\mathbb{P}_{\mu_\varepsilon}\right)_{\varepsilon>0} $. 
\noindent For $ T\in\mathbb N $ and a metric space $ (E,r) $ we define the time restriction $ R_T $ and time reversal operator $ r_T $:
\begin{align*}\label{timeresttimerev}
R_T:&C([0,\infty),E)\longrightarrow C([0,T],E),\omega\mapsto \omega_{|[0,T]}\\
r_T:&C([0,T],E) \longrightarrow C([0,T],E),\omega\mapsto \omega(T-\cdot).
\end{align*}
For a measure $ \mathbb P $ on $ C([0,\infty),E) $ we define $ \mathbb P^T:=\mathbb P\circ R_T^{-1}.$
%\begin{align*}\label{timeresttimerev}
%R_T:&C([0,\infty),\{\Phi_1<\infty\}\times\Rd{d})\longrightarrow C([0,T],\{\Phi_1<\infty\}\times\Rd{d}),\omega\mapsto \omega_{|[0,T]}\\
%r_T:&C([0,T],\{\Phi_1<\infty\}\times\Rd{d}) \longrightarrow C([0,T],\{\Phi_1<\infty\}\times\Rd{d}),\omega\mapsto \omega(T-\cdot).
%\end{align*}
%For a measure $ \mathbb P $ on $ C([0,\infty),\{\Phi_1<\infty\}\times\Rd{d}) $ we define $ \mathbb P^T:=\mathbb P\circ R_T^{-1}.$
We need two additional lemmata. Their proofs are elementary.
\begin{Le}\label{Timerestight}
	Let $ (E,r) $ be a metric space, $ \left(\mathbb P_n\right)_{n\in\mathbb N} $ be a family of Probability measures on $C\left([0,\infty),E\right) $ and $ \delta>0 $.
%	\begin{enumerate}
%		\item[(i)] If $ K\subseteq C\left([0,\infty),E\right) $ is totally bounded then $ R_T(K) $ is totally bounded in $ C\left([0,T],E\right) $ for all $ T\in\mathbb N $.
%		\item[(ii)] 
		If $ K_T\subseteq C\left([0,T],E\right)$ is a totally bounded set such that $ \inf_{n\in\mathbb N}\mathbb P_n^T(K_T)>1-\frac{\delta}{2^T} $ for all $T\in\mathbb N $. Then the set $ K=\bigcap_{T\in\mathbb N}R_T^{-1}K_T $ is totally bounded in $C\left([0,\infty),E\right) $ and it holds $ \inf_{n\in\mathbb N}\mathbb P_n(K)>1-\delta $.
\end{Le}

\begin{Le}\label{sumtight}
	Assume $ (E,\mathcal T) $ is a topological vector space, carrying the Borel $ \sigma $-algebra. Let $ X_n^i $, $i=1,2 $ be a $ E- $valued random variables on the probability space $ (\Omega_n,\mathcal F_n,\mathbb P_n) $, $ n\in\mathbb N$. Assume that the families $ \left(\mathbb P_n(X_n^i\in\cdot)\right)_{n\in\mathbb N} $, $ i=1,2 $, are tight. Then also the family $ \left(\mathbb P_n(X_n^1+X_n^2\in\cdot)\right)_{n\in\mathbb N} $ is tight. 
%	\begin{proof}
%		Let $ \delta>0 $. Choose $ K^i\subseteq E $, $ i=1,2 $, compact such that $ \inf_n P_n(X_n^i\in K^i)>1-\frac{\delta}{2}$. Observe that $ K^1+K^2=+(K^1\times K^2) $ is compact by Tychonoff theorem and $ \left\{X_n^1\in K^1\right\}\cap \left\{X_n^2\in K^2\right\}\subseteq\left\{X_n^1+X_n^2\in K^1+K^2\right\} $.
%		Hence it holds 
%		\[\inf_n P_n(X_n^1+X_n^2\in K^1+K^2)>1-\delta.\]
%	\end{proof}
\end{Le}

\begin{T}\label{Theoremtightness}
	Assume $ (\Phi_12),(\Phi_13),(\Phi_15)-(\Phi_19) $ and $ (\Phi_21)-(\Phi_27),(\Phi_29)-(\Phi_211) $.
	The family $ \left(\mathbb{P}_{\mu_\varepsilon}\right)_{\varepsilon>0} $ is tight as measures  on $C\left([0,\infty),\Rd{2d}\right) $.
	\begin{proof}
		In the following we always consider $\Rd{2d} $ to be equipped with the metric $ m $ from (\ref{metric}) and let $ T \in\mathbb N $ be arbitrary. 
		By Lemma \ref{Timerestight} it suffices to show that the family of time restrictions $ \left(\mathbb P^T_{\mu_\varepsilon}\right)_{\varepsilon>0} $ is tight for all $ T\in\mathbb N $. 
		For $ i\in\{1,...,d\} $ the functions $ f_i $, $ g_i $ from (\ref{fg}) induce measurable maps $ \hat f_i,\hat g_i$ defined by \[\hat{f}_i: C([0,T],\Rd{2d})\longrightarrow C([0,T],\Rd{}),\omega\mapsto f_i\circ\omega,\]
		analogous definition for $ \hat g_i $. 
%		Observe that the state space  
%		$ \{\Phi_1<\infty\}\times\Rd{d} $ is in general not complete.
		Due to the Arzel\`{a}-Ascoli theorem a set $A\subseteq C([0,T],\Rd{2d}) $ is totally bounded iff $\hat f_i(A),\hat g_i(A) \subseteq C([0,T],\Rd{}) $ are totally bounded for all $ i\in\{1,...,d\} $. Hence, it suffices to prove tightness separately for the following kind of measures on $ C([0,T],\Rd{}) $:
		\begin{equation}
		\textbf{1.}\left(\mathbb P^{T}_{\mu_\varepsilon} \circ \hat f_i^{-1}\right)_{\varepsilon>0},  i\in\{1,...,d\},\quad \textbf{2.} \left(\mathbb P^{T}_{\mu_\varepsilon} \circ \hat g_i^{-1}\right)_{\varepsilon>0}, i\in\{1,...,d\}.
		\end{equation}
		  In the following let  $i\in\{1,...,d\} $ and denote integration w.r.t. $ \mathbb P^{T}_{\mu_\varepsilon} $ by $ \mathbb E_\varepsilon^{T}$.
		  \noindent
		\begin{enumerate}[leftmargin=*,align=left]
			\item[\textbf{1.}] Consider the semimartingale decomposition from (\ref{Martingale2}):
			\begin{equation}\label{decfi}
				f_i(X_t,V_t)=M^{[f_i],\varepsilon}_t-\int_0^tL_\varepsilon f_i(X_r,V_r)\,dr+f_i(X_0,V_0),\quad t\in[0,T].
			\end{equation}
%			\begin{equation}\label{decfi}
%			f_i(X_t,V_t)=\underbrace{M^{[f_i],\varepsilon}_t}_{\alpha_t}-\underbrace{\int_0^t\partial_{x_i}\Phi_1(X_r,V_r)\,dr}_{\beta_t}+\underbrace{[f_i](X_0,V_0)}_{\gamma_t},\quad t\in[0,T].
%			\end{equation}
%As mentioned in Remark \ref{RemarkMartingale}\textit{(i)} one has to modify the random variable $ \int_0^TL_\varepsilon f_i(X_r,V_r)\,dr $ to obtain that each summand in (\ref{decfi}) is continuous in $ t $.
			This implies that $ \hat f_i $ can be written as the sum of the $ C([0,T],\Rd{}) $-valued random variables $\left (M^{[f_i],\varepsilon}_t\right )_{t\in[0,T]}  $, 
			$\left (\int_0^tL_\varepsilon f_i(X_r,V_r)\,dr\right )_{t\in[0,T]} $ and $ (f_i(X_0,V_0))_{t\in[0,T]}$, see also Remark \ref{RemarkMP}\textit{(i)}. Due to Lemma \ref{sumtight} it suffices to show separately that the laws of the single summands are tight. 
			We start with the family $ \mathbb P^{T}_{\mu_\varepsilon} \circ \left(\left(M^{[f_i],\varepsilon}_t\right)_{t\in[0,T]}\right)^{-1}$, $ \varepsilon>0 $. 
			Since the initial distributions of this family of measures are tight, it suffices to show a bound for the increments, see \cite{KS05}[Chapter 2, Problem 4.11].
			Therefore, let $ 0\leq s\leq t\leq T $. Since $ f_i^2 \in D(L_\varepsilon)$ and $ L_\varepsilon f_i \in L^4(\Rd{2d},\mu_{\varepsilon}) $, (\ref{QVP}) and (\ref{fi}) imply that the quadratic variation process of $ \left(M_t^{[f_i],\varepsilon}\right)_{t\in[0,T]}$ is given by a constant times $ t $.
			We obtain tightness by the following estimate which is due to the Burkholder-Davis-Gundy inequality,
			\begin{equation}
				\mathbb E_\varepsilon^{T}\left[(M_t^{[f_i],\varepsilon}-M_s^{[f_i],\varepsilon})^4\right]\leq %C\mathbb E_\varepsilon\left[\left(\int_s^t1(X_r,V_r)\,dr\right)^2\right]
				 C(t-s)^2.\label{Mf1}
			\end{equation}
%			For the variation part $ \mathbb P^{T}_{\mu_\varepsilon} \circ \left(\left(\int_0^tL_\varepsilon f_i(X_r,V_r)\,dr\right)_{t\in[0,T]}\right)^{-1}$, $ \varepsilon>0 $, we proceed similar.
			 Due to (\ref{fi}), the H\"older inequality and the fact that $  \mu_\varepsilon $ is invariant for $ \mathbb P_{\mu_\varepsilon} $ we find for the variation part $ \mathbb P^{T}_{\mu_\varepsilon} \circ \left(\left(\int_0^tL_\varepsilon f_i(X_r,V_r)\,dr\right)_{t\in[0,T]}\right)^{-1}$, $ \varepsilon>0 $, the following estimate implying tightness
			\begin{equation}
				\mathbb E_\varepsilon^{T}\left[\left(\int_s^tL_\varepsilon f_i(X_r,V_r)\,dr\right)^2\right]
%%				(t-s)\mathbb E_\varepsilon\left[\int_s^t\abs{\partial_{x_i}\Phi_1}^2(X_r,V_r)\,dr\right]\notag\\
%				=(t-s)\int_s^t \mathbb E_\varepsilon\left[\abs{\partial_{x_i}\Phi_1}^2(X_r,V_r)\right]\,dr\notag\\
%				=(t-s)\int_s^t \int_{\Rd{2d}}\abs{\partial_{x_i}\Phi_1}^2\,d\mu_{\varepsilon}\,dr\notag\\
				\leq (t-s)^2\mu_{\tilde{\Phi}_2}(\Rd{d})\int_{\Rd{d}}\abs{\partial_{i}\tilde \Phi_1}^2\,d\mu_{\tilde\Phi_1}.\label{Mf2}
			\end{equation}
			Tightness of the laws of the last summand follows by the weak convergence of the initial distributions and the continuity of $ f_i $.				
			We conclude that for  $ i\in\{1,...,d\} $ and $ T\in \mathbb N $
			 the family $ (\mathbb P^T_{\mu_\varepsilon}\circ \hat f_i^{-1})_{\varepsilon>0} $ is tight.
			\item[\textbf{2.}] It holds $ g_i\in D(L_\varepsilon)\cap D(\hat{L}_\varepsilon) $. Observe that  $ \mathbb P^T_{\mu_{\varepsilon}}\circ r_T^{-1} $ is associated with the adjoint semigroup $ \asg{\varepsilon}{,2} $, see \cite[Lemma 3.9(iii)]{GS15}, hence, $ \mathbb P^T_{\mu_{\varepsilon}}\circ r_T^{-1} =\hat{\mathbb P}^T_{\mu_{\varepsilon}}$. Explicit computation yields the following decomposition
%			, see also \cite[Equation (1)]{Tru2}, 
			\begin{align}\label{GLZ}
			g_i(X_t,V_t)-g_i(X_0,V_0)=&\frac{1}{2}\left(M_t^{g_i,\varepsilon}+\hat{M}_{T-t}^{g_i,\varepsilon}(r_T)-\hat M_T^{g_i,\varepsilon}(r_T)\right)\notag\\
			&+\frac{1}{2}\int_0^t(L_\varepsilon g_i-\hat L_\varepsilon g_i)(X_s,V_s)\,ds,\quad t\in[0,T].
			\end{align}
%			\begin{equation}\label{GLZ}
%				g_i(X_t,V_t)-g_i(X_t,V_0)=\frac{1}{2}\left(M_t^{g_i,\varepsilon}+\hat{M}_{T-t}^{g_i,\varepsilon}(\tau_T)-\hat M_T^{g_i,\varepsilon}(\tau_T)\right)+\frac{1}{2}\int_0^t(L_\varepsilon g_i-\hat L_\varepsilon g_i)(X_s,V_s)\,ds.
%			\end{equation}
			As above, we consider (\ref{GLZ}) as a decomposition of the random variable $ \hat g_i $. Tightness of $ \mathbb P^T_{\mu_{\varepsilon}}\circ\left((M_t^{g_i,\varepsilon})_{t\in[0,T]}\right)^{-1} $, $ \varepsilon>0 $, can be shown as in (\ref{Mf1}).
%			The first term on the right-hand side of (\ref{GLZ}) can also be estimated as in (\ref{Mf1}) since the quadratic variation process is the same as in (\ref{Mf1})
			For the summand $ \left(\hat{M}_{T-t}^{g_i,\varepsilon}(r_T)-\hat M_T^{g_i,\varepsilon}(r_T)\right)_{t\in[0,T]} $ we use $ \mathbb P^T_{\mu_{\varepsilon}}\circ r_T^{-1} =\hat{\mathbb P}^T_{\mu_{\varepsilon}}$. Since $ \left(\hat{M}_t^{g_i,\varepsilon}\right)_{t\in[0,T]} $ is a martingale w.r.t. $ \hat{\mathbb P}^T_{\mu_\varepsilon} $ tightness follows as (\ref{Mf1}).
%			:
%			\begin{equation*}
%				\mathbb E_\varepsilon^{T}\left[(\hat M_{T-s}^{g_i,\varepsilon}(r_T)-\hat M_{T-t}^{g_i,\varepsilon}(r_T))^4\right]
%				=\hat{\mathbb E}_\varepsilon\left[(\hat M_{T-s}^{g_i,\varepsilon}-\hat M_{T-t}^{g_i,\varepsilon})^4\right]\notag\\
%				\leq C(t-s)^2.
%			\end{equation*}
%			In the third line we applied the Burkholder-Davis-Gundy inequality. The quadratic variation process is determined by (\ref{QVP}) for $ \left(\hat L_\varepsilon ,D(\hat L_\varepsilon)\right)$ together with (\ref{gi2})
			Due to Proposition \ref{figi} we have for the last summand $\frac{1}{2}( L_\varepsilon g_i-\hat L_\varepsilon g_i)=-\partial_{x_i}\Phi_1 $, implying tightness of the laws $\mathbb P_{\mu_\varepsilon}^T\circ\left(\left( \int_0^t(L_\varepsilon g_i-\hat L_\varepsilon g_i)(Z_{s})\,ds\right)_{t\in[0,T]}\right)^{-1} $, $ \varepsilon>0 $, as in (\ref{Mf2}), which finishes the proof.
		\end{enumerate}
	\end{proof}
\end{T}
\noindent
%Let $ h_\varepsilon\in\mathcal H_\varepsilon $, $ \varepsilon>0 $, be a probability density w.r.t. $ \mu_{\varepsilon} $, such that $ \sup_{\varepsilon>0}\norm{h_\varepsilon}_{\mathcal H_\varepsilon}<\infty $.
%For $ \varepsilon>0 $ we define new probability measures $ \tilde{\mathbb P}_{h_\varepsilon\mu_{\varepsilon}}$ by $ \tilde{\mathbb P}_{h_\varepsilon\mu_{\varepsilon}}(A)=\mathbb E_\varepsilon\left[ 1_Ah_\varepsilon(X_0,V_0)\right]$, where $ A\in\mathcal B(C\left([0,\infty),\{\Phi_1<\infty\}\times\Rd{d}\right)) $ and $ \mathbb E_\varepsilon $ denotes the expectation w.r.t. $ \mathbb P_{\mu_\varepsilon} $ as in the previous proof. Observe that $ \tilde{\mathbb P}_{h_\varepsilon\mu_{\varepsilon}}$ coincides with the measure $\mathbb P_{h_\varepsilon\mu_{\varepsilon}}$ prescribed in the beginning of this section since it is also associated with the semigroup $ \sg{\varepsilon}{,2} $
%is the sense of (\ref{associated}).
%By Theorem \ref{Theoremtightness} the measures $ \mathbb P_{\mu_\varepsilon} $, $ \varepsilon>0 $, are tight. Now Hölders inequality implies tightness of $ \mathbb P_{h_\varepsilon\mu_{\varepsilon}}$, $\varepsilon>0 $.
 Combining Theorem \ref{fdd} and Theorem \ref{Theoremtightness} we obtain 
\begin{Co}\label{FinalConResult}
	Under the assumptions of Theorem \ref{fdd} and Theorem \ref{Theoremtightness} the measures $\left( \mathbb P_{h_\varepsilon\mu_{\varepsilon}}^X \right)_{\varepsilon>0 }$ on $ C\left([0,\infty),\Rd{d}\right) $ converge weakly to $ \mathbb P_{h_0\mu_0} $ for $ \varepsilon\to0 $.
	\begin{proof}
		By Theorem \ref{fdd} it suffices to prove tightness of $ \left(\mathbb P_{h_\varepsilon\mu_{\varepsilon}}^X\right)_{\varepsilon>0} $.
		The map $ P_X $ from (\ref{ProjX}) is continuous, hence, tightness of $ \left(\mathbb P_{h_\varepsilon\mu_{\varepsilon}}\right)_{\varepsilon>0} $ implies tightness of $ \left(\mathbb P_{h_\varepsilon\mu_{\varepsilon}}^X\right)_{\varepsilon>0} $. Now let $ \delta>0 $ and choose $ K\subseteq C\left([0,\infty),\Rd{2d}\right) $ compact s.t. $ \sup_{\varepsilon>0}\mathbb P_{\mu_\varepsilon}(K^c)\le \frac{\delta^2}{\sup_{\varepsilon>0}\norm{h_\varepsilon}^2_{L^2(\mu_\varepsilon)} }$. Again we denote by $ \mathbb E_\varepsilon $ integration w.r.t. $ \mathbb P_{\mu_\varepsilon} $.
		\begin{align*}
		\mathbb P_{h_\varepsilon\mu_{\varepsilon}}(K^c)&=\mathbb E_\varepsilon\left[1_{K^c}h_\varepsilon(X_0,V_0)\right]\leq \sqrt{\mathbb P_{\mu_\varepsilon}(K^c)}\norm{h_\varepsilon}_{L^2(\mu_\varepsilon)}\leq \delta.
		\end{align*}
	\end{proof}
\end{Co}

% 
% that $ \mathbb P_{\mu_{\varepsilon}}^X $ converge weakly to $ \mathbb P_{\mu_{0}}=\mathbb P_{\mu_{ \Phi_1}}$, i.e., the marginal distribution of the positions of the martingale solutions of the scaled Langevin equation converge weakly to the solution of the distorted Brownian motion.
%%
%i.e., the marginals of the positions of $ \mathbb P_{\mu_{\varepsilon^2}} $ and $ \mathbb P_{\mu_{\varepsilon^2}}\circ \left(\hat U_\varepsilon^{-1}\right)^{-1} $ are the same, which finishes our proof.

% !TeX spellcheck = en_US

\section[Overdamped limit of generalized stochastic Hamiltonian systems]{Overdamped limit of generalized stochastic Hamiltonian systems}\label{Overdampedlimit}
\setcounter{equation}{1}
Let us recall the scaled gsHs (\ref{sHseps1}), (\ref{sHseps2})
\begin{align*}
dX^\varepsilon_t&=\frac{1}{\varepsilon}\nabla\Phi_2(V^\varepsilon_t)dt,
%\tag{7.1a}
\\ dV^\varepsilon_t&=-\frac{1}{\varepsilon}\nabla\Phi_1(X^\varepsilon_t)dt-\frac{1}{\varepsilon^2}\nabla\Phi_2(V^\varepsilon_t)dt+\frac{1}{\varepsilon}\sqrt{2}dB_t,
%\tag{7.1a}
\end{align*}
%with generator given in (\ref{generatoroverdampedshs}) by
%
%%Observe that for $ \Phi_2(v)=\frac{1}{2}v^2 +\ln(\sqrt{2\pi}^d) $ we just recover (\ref{Langeps1}), (\ref{Langeps2}).
 We summarize our final result in the following theorem. To formulate the theorem define the map $  \tilde U_\varepsilon:\Rd{2d}\longrightarrow\Rd{2d},(x,v)\mapsto (x,\frac{v}{\varepsilon}), $ $ \varepsilon>0 $. In the following we denote by $ \mu $ the measure $ \mu_{\mathbf{\Phi}} $.
%Notation as in previous section 
%
%Let us describe how these results apply to the original problem, the convergence of the positions of the Langevin dynamics in the overdamped limit. 
%In the following we consider the case $ \Phi_2(v)=\frac{1}{2}v^2 +\ln(\sqrt{2\pi}^d) $. Observe that the assumption $ (\Phi_21)-(\Phi_211) $ hold for this particular choice of $ \Phi_2 $.
\begin{T}\label{finalresult}
Assume $ (\Phi_11)-(\Phi_19) $ and $ (\Phi_21)-(\Phi_211) $. Let $ \varepsilon>0 $, $ h_\varepsilon\in L^1(\Rd{2d},\mu)\cap L^2(\Rd{2d},\mu) $ and $ h\in L^1(\Rd{d},\mu_{\Phi_1})\cap L^2(\Rd{d},\mu_{\Phi_1}) $ be a probability densities w.r.t. $\mu $ and $ \mu_{\Phi_1} $, respectively. 
%Denote by $ \mathbb P_{h\mu_{\Phi_1}}$ the martingale solution from Corollary \ref{MPL} for the generator of the overdamped Langevin equation (\ref{Overdamped1}).
 Assume further that $  h_\varepsilon\mu $ converges weakly to $ h\mu_{\Phi_1}\otimes \delta_0 $ as $ \varepsilon\to0 $ and $ \sup_{\varepsilon>0}\int_{\Rd{2d}}h_\varepsilon^2d\mu<\infty $. There exists a weak solution $ \left(X^\varepsilon_t,V^\varepsilon_t\right)_{t\geq0} $ to (\ref{sHseps1}), (\ref{sHseps2}) with initial distribution $ h_\varepsilon\mu $. Furthermore, denote by $ \mathbb P_{h\mu_{\Phi_1}}$ the martingale solution to the generator of (\ref{Overdamped1}) from Corollary \ref{MPL}. Then the laws $ \mathcal L\left(\left(X^\varepsilon_t\right)_{t\geq0}\right) $, $ \varepsilon>0 $, converge weakly to $ \mathbb P_{h\mu_{\Phi_1}}$ as measures on $ C\left([0,\infty),\Rd{d}\right)$ as $ \varepsilon\to0 $.
\begin{proof}
Let $ \varepsilon>0 $ and recall $ \Phi_2^\varepsilon $, $ \left(L_\varepsilon,D(L_\varepsilon)\right) $, $ \sg{\varepsilon}{,2} $, $ \mu_\varepsilon $ and $ \mathcal H_\varepsilon $ from the beginning of Section \ref{CoFuAnOb}. 
The generator of (\ref{sHseps1}), (\ref{sHseps2}) is given by
\begin{equation}
L^\varepsilon_{\mathbf\Phi} f=\frac{1}{\varepsilon^2}\left(\Delta_vf-\nabla_v\Phi_2\cdot\nabla_vf\right)+\frac{1}{\varepsilon}\left(\nabla_v\Phi_2\cdot\nabla_xf-\nabla_x\Phi_1\cdot\nabla_vf\right), \quad f\in C^\infty_c(\{\Phi_1<\infty\}).
\end{equation}
We consider $\left( L^\varepsilon_{\mathbf\Phi},C^\infty_c(\{\Phi_1<\infty\})\right) $ as a linear operator on the space $\mathcal H=L^2\left(\Rd{2d},\mu\right) $.
 Define the unitary transformation $ U_\varepsilon : \mathcal H\longrightarrow\mathcal H_{\varepsilon} , f\mapsto f\circ \tilde U_\varepsilon$. The map $ U_\varepsilon $ and the adjoint $ U_\varepsilon^* $ leave the set $ C^\infty_c(\{\Phi_1<\infty\}) $ invariant. Furthermore, we obtain the unitary equivalence
\begin{equation}
\left( U_\varepsilon^* L^1_{\mathbf{\Phi}^\varepsilon} U_\varepsilon,C^\infty_c(\{\Phi_1<\infty\})\right)=\left(L^{\varepsilon}_{\mathbf\Phi},C^\infty_c(\{\Phi_1<\infty\})\right).
\end{equation}
By Lemma \ref{Lpgen} an extension of $\left(L^1_{\mathbf\Phi_\varepsilon},C^\infty_c(\{\Phi_1<\infty\})\right) $ is the generator of the semigroup $ \left( T_{t,2}^{\varepsilon}\right)_{t\geq0}$.
Hence, due to \cite{G85}[Chapter 2, Lemma 3.17] an extension of the operator $\left( L^\varepsilon_{\mathbf\Phi},C^\infty_c(\{\Phi_1<\infty\})\right) $ is the generator of the sub-Markovian strongly continuous contraction semigroup on $ \mathcal H $ given by  $\left(S_t^\varepsilon\right)_{t\geq0}=\left( U_\varepsilon^* T_{t,2}^{\varepsilon} U_\varepsilon\right)_{t\geq0} $.
Define further
\begin{equation*}
\hat U_\varepsilon:C([0,\infty),\Rd{2d})\longrightarrow C([0,\infty),\Rd{2d}), (x_t,v_t)_{t\geq 0}\mapsto \left(\tilde U_\varepsilon(x_t,v_t)\right)_{t\geq 0}.
\end{equation*} 
Observe that $U_\varepsilon h_\varepsilon $ is a probability density w.r.t.~$ \mu_\varepsilon $. Let $ \mathbb P_{(U_\varepsilon h_\varepsilon)\mu_\varepsilon} $
be the martingale solution to $ \left(L^1_{\mathbf\Phi^\varepsilon,2},D\left(L^1_{\mathbf\Phi^\varepsilon,2}\right)\right) $ with initial distribution $ (U_\varepsilon h_\varepsilon)\mu_\varepsilon $ from the last section. One easily checks that the measure $\tilde{\mathbb{P}}_{h_\varepsilon\mu}:= \mathbb P_{(U_\varepsilon h_\varepsilon)\mu_{\varepsilon}}\circ \left(\hat U_\varepsilon
\right)^{-1} $ has initial distribution given by $ h_\varepsilon\mu $ and is associated with the sub-Markovian semigroup $ \left(S_t^{\varepsilon}\right)_{t\geq0} $ in the sense of (\ref{associated}). Hence, due to \cite[Lemma 5.1]{CG10} the measure $\tilde{\mathbb{P}}_{h_\varepsilon\mu}$ is a martingale solution to the generator of $ \left(S_t^{\varepsilon}\right)_{t\geq0} $. Furthermore, one can argue as in Remark \ref{weaksolution} to obtain weak solutions $ (X_t^\varepsilon,V_t^\varepsilon)_{t \geq 0} $ from $\tilde{\mathbb{P}}_{h_\varepsilon\mu}$ such that for the law of $ (X_t^\varepsilon)_{t \geq 0} $ it holds $ \mathcal L\left(\left(X^\varepsilon_t\right)_{t\geq0}\right)=\tilde{\mathbb P}_{h_\varepsilon\mu}\circ P_X^{-1}$.
Observe that $ \tilde{\mathbb P}_{h_\varepsilon\mu}\circ P_X^{-1} =\mathbb P_{(U_\varepsilon h_\varepsilon)\mu_{\varepsilon}}\circ P_X^{-1}$. To apply Corollary \ref{FinalConResult} we have to guarantee that the assumptions of Theorem \ref{fdd} are fulfilled, i.e., we have show that $(U_\varepsilon h_\varepsilon)\mu_\varepsilon $, $ \varepsilon>0 $, converges weakly to $ h\mu_{\Phi_1}\otimes \delta_0 $ as $ \varepsilon\to0 $.
Let $ f:\Rd{2d}\longrightarrow \Rd{} $ be continuous and bounded. Observe that the functions $ g_\varepsilon $ defined by $ g_\varepsilon(x,v)=f(x,\varepsilon v) $ converge uniformly on compact sets to the function $ g(x,v)=f(x,0) $, $ (x,v)\in\Rd{2d} $. Hence, by the transformation formula we obtain
\begin{equation*}
\int_{\Rd{2d}}f (U_\varepsilon h_\varepsilon)d\mu_\varepsilon=\int_{\Rd{2d}}g_\varepsilon h_\varepsilon d\mu
=\int_{\Rd{2d}}(g_\varepsilon-g) h_\varepsilon d\mu+\int_{\Rd{2d}}g h_\varepsilon d\mu.
\end{equation*}
It suffices to prove that the first term in the last expression converges to zero as $ \varepsilon\to 0 $.  By assumption the measures $ h_\varepsilon\mu $, $ \varepsilon>0 $ converge weakly, in particular, they are tight. Hence by the boundedness of $ f $ and the considerations above we conclude 
\begin{equation*}
\int_{\Rd{2d}}f(U_\varepsilon h_\varepsilon)d\mu_\varepsilon\xrightarrow{\varepsilon\to 0}\int_{\Rd{2d}}fd h\mu_{\Phi_1}\otimes \delta_0.
\end{equation*}
Hence, we can apply Corollary \ref{FinalConResult} and conclude that $ 
\tilde{\mathbb P}_{h_\varepsilon\mu}\circ P_X^{-1} =\mathbb P_{(U_\varepsilon h_\varepsilon)\mu_{\varepsilon}}\circ P_X^{-1} $ converge weakly to $ \mathbb P_{h\mu_{\Phi_1}} $ which finishes the proof.
\end{proof} 
\end{T}
\begin{Rem}
	Recall the objects $ \mathcal H $, $ U_\varepsilon^* $, $ \left(S_t^\varepsilon\right)_{t\geq0} $, $ \varepsilon>0 $, from the previous proof.
	Via the maps $ \Psi_\varepsilon $ from (\ref{Psi}) one directly obtains $ \mathcal H\xrightarrow{\left(\Gamma_\varepsilon\right)_{\varepsilon>0}}\mathcal H_{\Phi_1} $, where $ \Gamma_\varepsilon:\mathcal C\longrightarrow\mathcal H, f\mapsto U_\varepsilon^*\circ \Psi_\varepsilon(f) $. Furthermore, we obtain that the semigroups $ \left(S_t^\varepsilon\right)_{t\geq0} $ converge to $ \sg{\Phi_1}{} $ along $ \mathcal H\xrightarrow{\left(\Gamma_\varepsilon\right)_{\varepsilon>0}}\mathcal H_{\Phi_1} $. This follows directly from the fact that the properties (\ref{normconv}), (\ref{weakconv}) are preserved by the unitary map $ U_\varepsilon^* $.
\end{Rem}

% !TeX spellcheck = en_US
\section*{Acknowledgement}
The second author thanks the department of Mathematics at the University of Kaiserslautern for financial support in the form of a fellowship. 

%\printbibliography
\bibliography{references} 	
\end{document}